\pdfoutput=1 

\documentclass[a4paper, 10pt]{amsart}

\usepackage{amssymb, amsmath, bbm, dsfont, mathrsfs, enumitem, url}
\usepackage{wasysym, mathtools, subfigure}
\usepackage[dvipsnames, usenames]{color}
\usepackage[colorlinks=true,
        raiselinks=true,
        pagecolor=RawSienna,
        linkcolor=MidnightBlue,
        citecolor=Mahogany,
        urlcolor=ForestGreen,
        pdfauthor=D. Chatterjee,
        pdftitle={Stochastic receding horizon control with bounded control inputs: a vector-space approach},
        pdfkeywords={stochastic MPC, unbounded noise, bounded/saturated control inputs},
        pdfsubject={Technical Report},
        plainpages=false]{hyperref}
\usepackage[pdftex]{graphicx}
\usepackage{sidecap, wrapfig}

\usepackage{pxfonts}
\definecolor{mred}{rgb}{0.6, 0, 0}
\definecolor{mgreen}{rgb}{0, 0.5, 0}
\definecolor{mblue}{rgb}{0, 0, 0.5}
\definecolor{mcyan}{rgb}{0, 0.5, 0.5}

\let\oldsqrt\sqrt
\def\sqrt{\mathpalette\DHLhksqrt}
\def\DHLhksqrt#1#2{%
\setbox0=\hbox{$#1\oldsqrt{#2\,}$}\dimen0=\ht0
\advance\dimen0-0.2\ht0
\setbox2=\hbox{\vrule height\ht0 depth -\dimen0}%
{\box0\lower0.4pt\box2}}

\makeatletter
\newcommand\@dotsep{4.5}
\def\@tocline#1#2#3#4#5#6#7{\relax
  \ifnum #1>\c@tocdepth 
  \else
    \par \addpenalty\@secpenalty\addvspace{#2}%
    \begingroup \hyphenpenalty\@M
    \@ifempty{#4}{%
      \@tempdima\csname r@tocindent\number#1\endcsname\relax
    }{%
      \@tempdima#4\relax
    }%
    \parindent\z@ \leftskip#3\relax \advance\leftskip\@tempdima\relax
    \rightskip\@pnumwidth plus1em \parfillskip-\@pnumwidth
    #5\leavevmode\hskip-\@tempdima #6\relax
    \leaders\hbox{$\m@th
      \mkern \@dotsep mu\hbox{.}\mkern \@dotsep mu$}\hfill
    \hbox to\@pnumwidth{\@tocpagenum{#7}}\par
    \nobreak
    \endgroup
  \fi}
\makeatother

\DeclareMathOperator{\trace}{\mathbf{tr}}
\DeclareMathOperator{\sinc}{sinc}
\DeclareMathOperator{\linspan}{span}
\DeclareMathOperator{\Gammaf}{Gamma}

\DeclareMathOperator{\erf}{erf}
\DeclareMathOperator{\erfc}{erfc}

\DeclareMathOperator{\diag}{diag}

\DeclareMathOperator{\sgn}{sgn}

\newcommand{\Let}{\coloneqq}
\newcommand{\teL}{\eqqcolon}

\def\NN{\mathbb{N}}
\def\Nz{\mathbb{N}_0}
\def\R{\mathbb{R}}
\def\C{\mathbb{C}}
\def\posR{\mathbb{R}_{\ge 0}}
\def\EE{\mathbb{E}}
\def\PP{\mathbb{P}}

\def\U{\mathbb{U}}
\def\W{\mathbb{W}}
\def\Z{\mathbb{Z}}

\def\x0{x_0}
\def\xz{x_0}

\def\ii{\mathrm i}
\newcommand{\HH}{\mathcal H}
\newcommand{\iss}{\textsc{iss}}
\newcommand{\ClassKL}{\mathcal{KL}}
\newcommand{\ClassKinfty}{\mathcal{K}_\infty}
\newcommand{\fa}{\forall\,}
\newcommand{\wh}{\widehat}
\newcommand{\mcal}{\mathcal}

\newcommand{\ee}{{\mathfrak e}}

\renewcommand{\subset}{\subseteq}
\renewcommand{\le}{\leqslant}
\renewcommand{\ge}{\geqslant}
\newcommand{\nn}{\nonumber}
\newcommand{\mc}{\mathcal}

\newcommand{\wt}{\widetilde}
\newcommand{\transp}{^\mathsf{T}}
\newcommand{\mrm}{\mathrm}

\newcommand{\ExEnd}{\hspace{\stretch{1}}{$\triangle$}}
\newcommand{\DefEnd}{\hspace{\stretch{1}}{$\Diamond$}}
\newcommand{\AssmpEnd}{\hspace{\stretch{1}}{$\diamondsuit$}}

\newcommand{\smat}[1]{\left[\begin{matrix} #1 \end{matrix}\right]}

\newcommand{\tr}[1]{\mathbf{tr}\!\left(#1\right)}
\newcommand{\inprod}[2]{\left\langle{#1},{#2}\right\rangle}
\newcommand{\lambdamax}[1]{\lambda_{\text{max}}(#1)}

\newcommand{\norm}[1]{\left\lVert{#1}\right\rVert}

\newcommand{\epower}[1]{\mathrm e^{#1}}
\newcommand{\eps}{\varepsilon}
\newcommand{\indic}[1]{\mathbf 1_{#1}}
\newcommand{\secref}[1]{\S\ref{#1}}

\numberwithin{equation}{section}
\newtheoremstyle{dcstyle}{10pt}{10pt}{\slshape}{}{\bfseries}{.}{ }{}
\newtheoremstyle{nonum}{10pt}{10pt}{}{}{\itshape}{.}{ }{\thmname{#1}\thmnote{ (\mdseries #3)}}
\theoremstyle{dcstyle}
\newtheorem{theorem}{Theorem}[section]
\newtheorem{corollary}[theorem]{Corollary}

\newtheorem{proposition}[theorem]{Proposition}

\theoremstyle{definition}
\newtheorem{defn}[theorem]{Definition}
\newtheorem{assumption}[theorem]{Assumption}

\theoremstyle{remark}

\newtheorem{example}[theorem]{Example}

\theoremstyle{nonum}

\newtheoremstyle{nonumF}{}{}{\itshape}{}{\bfseries}{.}{ }{\thmname{#1}\thmnote{ (\mdseries #3)}}
\theoremstyle{nonumF}
\newtheorem{fact}{Facts about Special Functions}

\allowdisplaybreaks

\title[Stochastic receding horizon control with bounded control inputs]{Stochastic receding horizon control with bounded control inputs: a vector space approach}
\thanks{This research was partially supported by the Swiss National Science Foundation under grant 200021-122072, and the FeedNetBack project FP7-ICT-223866 (\url{www.feednetback.eu}).}
\author[D.~Chatterjee]{Debasish Chatterjee}
\author[P.~Hokayem]{Peter Hokayem}
\author[J.~Lygeros]{John Lygeros}
\address{Automatic Control Laboratory\\Physikstrasse 3\\ETH Z\"urich\\8092 Z\"urich\\Switzerland}
\email{\{chatterjee,hokayem,lygeros\}@control.ee.ethz.ch}

\begin{document}


\begin{abstract}
We design receding horizon control strategies for stochastic discrete-time linear systems with additive (possibly) unbounded disturbances, while obeying hard bounds on the control inputs. We pose the problem of selecting an appropriate optimal controller on vector spaces of functions and show that the resulting optimization problem has a tractable convex solution. Under the assumption that the zero-input and zero-noise system is asymptotically stable, we show that the variance of the state is bounded when enforcing hard bounds on the control inputs, for any receding horizon implementation.  Throughout the article we provide several examples that illustrate how quantities needed in the formulation of the resulting optimization problems can be calculated off-line, as well as comparative examples that illustrate the effectiveness of our control strategies.
\end{abstract}

\maketitle

\begin{minipage}[c]{0.85\textwidth}
\begin{small}
	\tableofcontents
\end{small}
\end{minipage}

\section{Introduction}
\label{s:intro}

Receding horizon control is a popular paradigm for designing control policies. In the context of deterministic systems it has received a considerable amount of attention over the last two decades, and significant advancements have been made in terms of its theoretical foundations as well as industrial applications. The motivation comes primarily from the fact that receding horizon control yields tractabile control laws for deterministic systems in the presence of constraints, and has consequently become popular in the industry. The counterpart in the context of stochastic systems, however, is a relatively recent development. In this article we solve the problem of stochastic receding horizon control for linear systems subject to additive (possibly) unbounded disturbances and hard norm bounds on the control inputs, over a class of feedback policies. Methods for guaranteeing hard bounds on the control inputs, within our context, while ensuring tractability of the underlying optimization problem are, to our knowledge, not available in the current literature. Preliminary results in this direction were reported in \cite{ref:petercdc09}.

In the deterministic setting, the receding horizon control scheme is dominated by worst-case analysis relying on robust control and robust optimization methods, see, for example, \cite{ref:bertsekas05survey, MayneRawlingsRaoScokaert-00,BemporadMorari-99, ref:lazarbemporad07, ref:maciejowskibk, ref:blanchini1999sic, FukushimaBitmead-05, ref:JunYanBitmead, RichardsHow-05} and the references therein. The central idea is to synthesize a controller based on the bounds of the noise such that a certain target set becomes invariant with respect to the closed-loop dynamics. However, such an approach tends to yield rather conservative controllers and large infeasibility regions. Moreover, assigning an a priori bound to the noise seems to demand considerable insight. A stochastic model of the noise is a natural alternative approach to this problem: the conservativeness of worst-case controllers may be reduced, and one may not need to impose any a priori bounds on the maximum magnitude of the noise. In~\cite{BertsimasBrown-07}, the authors reformulate the stochastic programming problem as a deterministic one with bounded noise and solve a robust optimization problem over a finite horizon, followed by estimating the performance when the noise is unbounded but takes high values with low probability (as in the Gaussian case). In~\cite{PrimbsSung-09} a slightly different problem is addressed in which the noise enters in a multiplicative manner, and hard constraints on the states and control inputs are relaxed to constraints resembling the integrated chance constraints of \cite{ref:haneveldICC} or risk measures in mathematical finance. Similar relaxations of hard constraints to soft probabilistic ones have also appeared in~\cite{ref:CannonKouvaritakisWu-08} for both multiplicative and additive noise inputs, as well as in~\cite{OldewurtelJonesMorari-08} for additive noise inputs. There are also other approaches, e.g., those employing randomized algorithms as in~\cite{BlackmoreWilliams-07, batinaPhDthesis, MaciejowskiLecchiniLygeros-05}. Related lines of research can be found in~\cite{vanHessemFullSolution, vanHessem2006} dealing with constrained model predictive control (MPC) for stochastic systems motivated by industrial applications, in~\cite{ref:ramponinstable, ref:batina02, ref:Stoorvogel} dealing with stochastic stability, in \cite{SkafBoyd-Qdesign-09} dealing with Q-design, in, e.g., \cite{lavretsky2007stable, LavretskyHovakimyanCalise-03} dealing with alternative approaches to control under actuator constraints and neural-network approximation. The articles \cite{AgarwalCinquemaniChatterjeeLygeros-09, ref:cinquemaniagarwal} deal with a formulation that allows probabilistic state constraints but not hard input constraints, and is hence complementary to the approach in the present article, and \cite{ref:acc10} treats the case of output feedback. . Finally, note that probabilistic constraints on the controllers naturally raise difficult questions on what actions to take when such constraints are violated, see~\cite{ref:recstrat} and~\cite{ref:excursion} for partial solutions to these issues.

The main contributions of the article are as follows: We give a tractable, convex, and globally feasible solution to the finite-horizon stochastic linear quadratic (LQ) problem for linear systems with possibly unbounded additive noise and hard constraints on the elements of the control policy. Within this framework one has two directions to pursue in terms of controller design, namely, a posteriori bounding the standard LQG controller, or employing certainty-equivalent receding horizon controller. While the former direction explicitly incorporates some aspects of feedback, the synthesis of the latter involves control constraints and implicitly incorporates the notion of feedback. Our choice of feedback policies explores the middle ground between these two choices: we explicitly incorporate both the control bounds and feedback at the design phase. More specifically, we adopt a policy that is affine in certain bounded functions of the past noise inputs. The optimal control problem is lifted onto general vector spaces of candidate control functions from which the controller can be selected algorithmically by solving a convex optimization problem. Our novel approach does not require artificially relaxing the hard constraints on the control input to soft probabilistic ones (to ensure large feasible sets), and still provides a globally feasible solution to the problem. Minimal assumptions of the noise sequence being i.i.d and having finite second moment are imposed. The effect of the noise appears in the convex optimization problem as certain fixed cross-covariance matrices, which may be computed off-line and stored.

Once tractability of the optimization problem is ensured, we employ the resulting control policy in a receding horizon scheme. Under our policies the closed-loop system is in general not necessarily Markovian, and as a result stability of the closed-loop system is not immediate. In fact, we can no longer appeal directly to standard Foster-Lyapunov methods. We establish that our receding horizon control scheme provides stability under the assumption that the zero-input and zero-noise system is asymptotically stable. We provide examples that demonstrate the effectiveness of our policies with respect to standard methods such as certainty-equivalent MPC, standard unconstrained LQG and saturated LQG control. These examples show that our policies perform no worse than the standard unconstrained LQG controller in the absence of control constraints, and outperform the certainty-equivalent MPC as well as the saturated LQG control by a significant margin.

Our mechanism for selection of a policy consists of two steps: The first concerns the structure of our policies, and is motivated by preceding work in robust optimization and MPC \cite{ref:loefberg03, ref:ben-tal04, ref:goulart06}. The second concerns the procedure for selection of an optimal policy from a general vector space of candidate control functions, and is inspired by approximate dynamic programming techniques \cite{ref:bertsekasNDP, ref:rantzerRelaxingDP, ref:schweitzerPolyApprox, ref:vanRoyLPtoDP, ref:powellADP}. With respect to the first step, our policies are more general compared to those in \cite{ref:loefberg03, ref:ben-tal04, ref:goulart06}. With respect to the second, the selection procedure of our policies consists of a one-step tractable static optimization program.

The rest of this article is organized as follows. In Section \ref{s:ps} we state the main problem to be solved in the most general form. In Section \ref{s:mainres} we provide a tractable solution to the finite horizon optimization problem on general vector spaces. This result is specialized to various classes of noise and input constraint sets in Section \ref{s:if}. Stability of receding horizon implementations of the obtained closed-loop policy is shown in Section \ref{sec:stability}, and input-to-state stability properties are discussed in Section \ref{sec:iss}. We provide a host of numerical examples that illustrate the effectiveness of our approach in Section \ref{sec:examples}. Finally, we conclude in Section \ref{sec:conclusions} with a discussion on future research directions.


\subsection*{Notation}
Hereafter, $\NN \Let \{1,2,\ldots\}$ is the set of natural numbers, $\Nz \Let \NN \cup \{0\}$, $\Z$ is the set of all the integers, $\posR$ is the set of nonnegative real numbers, and $\C$ denotes the set of complex numbers. We let $\indic{A}(\cdot)$ denote the indicator function of a set $A$, and $\mathbf I_{n\times n}$ and $\mathbf 0_{n\times m}$ denote the $n$-dimensional identity matrix and $n\times m$-dimensional zeros matrix, respectively. Let $\norm{\cdot}$ denote the standard Euclidean norm, and $\norm{\cdot}_p$ denote the usual $\ell_p$ norms. Also, let $\EE_{\xz}[\cdot]$ denote the expected value given $\xz$, and $\tr{\cdot}$ denote the trace of a matrix. If $M_1$ and $M_2$ are two matrices with the same number of rows, we employ the standard notation $[M_1\mid M_2]$ for the matrix obtained by stacking the columns of $M_1$ followed by the columns of $M_2$. For a given symmetric $n$-dimensional matrix $M$ with real entries, let $\{\lambda_i(M)\mid i=1, \ldots, n\}$ be the set of eigenvalues of $M$, and let $\lambda_{\rm max}(M) \Let \max_i\lambda_i(M)$ and $\lambda_{\text{min}}(M) \Let \min_i\lambda_i(M)$. Finally, for a random vector $X$ let $\Sigma_X$ denote the matrix $\EE\bigl[XX\transp\bigr]$ and $\mu_X$ denote the vector $\EE\bigl[X\bigr]$.


\section{Problem Statement}
\label{s:ps}

%

Consider the following discrete-time stochastic dynamical system:
\begin{equation}
\label{eq:system}
	x_{t+1} = \bar A x_t + \bar B u_t + w_t, \qquad t\in\NN_0,
\end{equation}
where $x_t\in\R^{n}$ is the state, $u_t$ is the control input taking values in some given control set $\bar \U\subset\R^m$ to be defined later, $\bar A\in\R^{n\times n}$, $\bar B\in\R^{n\times m}$, and $(w_t)_{t\in\Nz}$ is a sequence of stochastic noise vectors with $w_t\in\W\subset\R^n$. We assume that the initial condition $x_0$ is given and that, at any time $t$, $x_t$ is observed perfectly. We do not assume that the components of the noise $w_t$ are uncorrelated, nor that they have zero mean; this effectively means that $w_t$ may be of the form $\bar F w_t' + b$ for some noise $w_t'\in\R^p$ whose components are uncorrelated or mutually independent, $F\in\R^{n\times p}$, and $b\in\R^n$. Without loss of generality we shall stick to the simpler notation of~\eqref{eq:system} throughout this article. The results readily extend to the general case of $w_t =\bar  F w_t' + b$, as can be seen in~\cite{ref:petercdc09}.

Generally, a  \emph{control policy} $\pi$ is a sequence $(\pi_0, \pi_1, \pi_2, \ldots)$ of Borel measurable maps $\pi_t:\underbrace{\R^n\times\cdots\times\R^n}_{k(t)-\text{ times}}\to\bar \U,\;t\in\Nz$. Policies of finite length such as $(\pi_t, \pi_{t+1}, \ldots, \pi_{t+N-1})$ will be denoted in the sequel by $\pi_{t:t+N-1}$. 

Fix an optimization horizon $N\in\NN$ and let us consider the following objective function at time $t$ given the state $x_t$:
\begin{equation}
\label{e:objfn}
	V_t \Let \EE\Biggl[\sum_{k=0}^{N-1}\bigl(x_{t+k}\transp Q_{k} x_{t+k}+u_{t+k}\transp R_{k} u_{t+k}\bigr) + x_{t+N}\transp Q_Nx_{t+N}\,\Bigg|\,x_t\Biggr],
\end{equation}
where $Q_t > 0, R_t > 0, Q_N > 0$ are some given symmetric matrices of appropriate dimension.
At each time instant $t$, we are interested in minimizing \eqref{e:objfn} over the class of causal state feedback strategies $\Pi$ defined as:
\begin{equation}
\label{e:policies}
\smat{u_t \\ u_{t+1}\\ \vdots \\ u_{t+N-1}} = \left[ \begin{array}{l} \pi_t(x_t) \\
\pi_{t+1}(x_t,x_{t+1}) \\ \vdots \\ \pi_{t+N-1}(x_t,x_{t+1},\cdots,x_{t+N-1})\end{array}\right],
\end{equation}
for some measurable functions $\pi_{t:t+N-1}\Let\{\pi_t,\cdots, \pi_{t+N-1}\}\in\Pi$, while satisfying $u_t \in \bar\U$ for each $t$. The \emph{receding horizon control} procedure for a given control horizon $N_c\in\{1, \ldots, N\}$ and time $t$ can be described as follows:
\begin{itemize}[label=(\alph*), leftmargin=*, align=right]
    \item[(a)]   measure the state $x_t$;

    \item[(b)] determine an admissible optimal feedback control policy, say $\pi^*_{t:t+N-1}\in\Pi$, that minimizes the $N$-stage cost function \eqref{e:objfn} starting from time $t$, given the measured initial condition $x_t$;
    \item[(c)] apply the first $N_c$ elements $\pi^*_{t:t+N_c-1}$ of the policy $\pi^*_{t:t+N-1}$;
	\item[(d)] increase $t$ to $t+N_c$, and go back to step (a).
\end{itemize}
In this context, if $N_c = 1$ then this is usual MPC, and if $N_c = N$, then it is usually known as rolling horizon control.

Since both the system \eqref{eq:system} and cost \eqref{e:objfn} are time-invariant, it is enough to consider the problem of minimizing the cost for $t = 0$. 
In view of the above we consider the problem:
\begin{equation}
\label{eq:problem}
\begin{aligned}
\min_{\pi_{0:N-1}\in\Pi} \bigl\{V_0\,\big|\,\text{dynamics \eqref{eq:system}},  \text{ and } u_t\in\bar \U \text{ for each } t\bigr\}.
\end{aligned}
\end{equation}
If feasible, the problem \eqref{eq:problem} generates an optimal sequence of feedback control laws $\pi^*=\left\{\pi^*_0,\cdots, \pi^*_{N-1}\right\}$.

The evolution of the system (\ref{eq:system}) over a single optimization horizon $N$ can be described in a compact form as follows:
\begin{equation}
\label{eq:compactdyn}
x=A\x0+Bu+Dw,
\end{equation}
where
$$x\Let \begin{bmatrix}
x_0 \\ x_1  \\ \vdots \\ x_N
\end{bmatrix},\qquad u\Let \begin{bmatrix}
u_0 \\ u_1 \\ \vdots \\ u_{N-1}
\end{bmatrix},\qquad
w \Let
\begin{bmatrix}
w_0 \\ w_1 \\ \vdots \\ w_{N-1}
\end{bmatrix},\qquad A \Let
\begin{bmatrix}
\mathbf I_{n\times n} \\ \bar A \\ \vdots \\\bar  A^N
\end{bmatrix},$$
$$
B \Let
\begin{bmatrix}
\mathbf 0_{n\times m} &\cdots &\cdots & \mathbf 0_{n\times m} \\
\bar B &\ddots &&\vdots\\
\bar A\bar B &\bar  B &\ddots & \vdots\\
\vdots && \ddots &\mathbf 0_{n\times m}\\
\bar A^{N-1}\bar  B & \cdots & \bar A\bar B&\bar  B
\end{bmatrix}, \qquad
D \Let
\begin{bmatrix}
\mathbf 0_{n\times n} &  \cdots &\cdots& \mathbf 0_{n\times n} \\
\mathbf I_{n\times n} & \ddots & &\vdots \\
\bar A & \mathbf{I}_{n\times n} &\ddots &\vdots \\
\vdots && \ddots & \mathbf 0_{n\times n} \\
\bar A^{N-1} &\cdots & \bar A& \mathbf I_{n\times n}
\end{bmatrix}.
$$
Using the compact notation above, the optimization Problem~\eqref{eq:problem} can be rewritten as follows:
\begin{equation}
\label{eq:problem1}
\begin{aligned}
\min_{\pi_{0:N-1}\in\Pi}   \quad&\bigl\{\EE_{\xz}\bigl[ x\transp  Q x+ u\transp  R u\bigr]\,\big|\, \text{dynamics \eqref{eq:compactdyn}}, u\in\U\bigr\},
\end{aligned}
\end{equation}
where $Q = \diag\{Q_0, \ldots, Q_{N}\}$, $R = \diag\{R_0,\ldots, R_{N-1}\}$, and $\U\Let\underbrace{\bar \U\times\ldots\times \bar \U}_{N-\text{times}}$.


\section{Main Result}
\label{s:mainres}
	We require that our controller is selected from a vector space of candidate controllers spanned by a given set of ``simple'' basis functions. The precise algorithmic selection procedure is based on the solution to an optimization problem. The basis functions may represent particular types of control functions that are easy or inexpensive to implement, e.g., minimum attention control \cite{ref:brockettMinAttention}, or may be the only ones available for a specific application. For instance, piecewise constant policy elements with finitely many elements in their range may be viewed as controllers that can provide only finitely many values; this may be viewed as an extended version of a bang-bang controller, or as a hybrid controller with a finite control alphabet.

	More formally, let $\HH$ be a nonempty separable vector space of functions with the control set $\U$ as their range, i.e., $\HH$ is the linear span of measurable functions $\ee^\nu:\W\to\U$, where $\nu\in\mathcal I$ - an ordered countable index set (see \cite{Luenberger-69} for more details). As mentioned above, the elements of $\HH$ may be linear combinations of typical ``simple'' controller functions for $t = 0, 1, \ldots, N-1$. We are interested in policies of the form $u_t = \eta_t + \sum_{i=0}^{t-1}\psi_{t, i}(w_i)$, where $\eta_t$ is an $m$-dimensional vector and each component of the $m$-dimensional vector-valued function $\psi_{t, i}$ is a member of $\HH$. Although this feedback function is directly from the noise, since the state is assumed to be perfectly measured, from the system dynamics~\eqref{eq:system} it follows at once that this controller $u_t$ is actually a feedback from all the states $x_0, \ldots, x_t$. Indeed, in the spirit of~\cite{ref:loefberg03, ref:ben-tal04, ref:goulart06, SkafBoyd-Affine-09} we have
	\begin{align*}
		u_0 & = \eta_0,\\
		u_1 & = \eta_1 + \psi_{1, 0}(x_1 -\bar  Ax_0 -\bar  B\eta_0),\\
		u_2 & = \eta_2 + \psi_{2, 0}(x_1 -\bar  Ax_0 -\bar  B\eta_0) + \psi_{2, 1}\bigl(x_2 -\bar  Ax_1 - \bar B\bigl(\eta_1 + \psi_{1, 0}(x_1 - \bar  A x_0 - \bar  B \eta_0)\bigr)\bigr),\\
		& \vdots
	\end{align*}
	In other words, by construction, $u_t$ is  generally a nonlinear feedback controller depending on the past $t$ states.\footnote{Note that the controller input at time $t$ is non-Markovian  as it is a function of the state vectors at all the previous times and not just on $x_{t-1}$.} Also by construction, it is causal.

	Our general control policy can now be expressed as the vector
	\begin{equation}
	\label{e:genpolicy}
		 u  = \eta + \varphi( w ) \Let
		\begin{bmatrix}
			\eta_0\\
			\eta_1\\
			\vdots\\
			\eta_{N-1}
		\end{bmatrix} +
		\begin{bmatrix}
			\varphi_0\\
			\varphi_1(w_0)\\
			\vdots\\
			\varphi_{N-1}(w_0, w_1, \ldots, w_{N-2})
		\end{bmatrix},
	\end{equation}
	where,
	\begin{itemize}[label=\textbullet, leftmargin=*]
		\item $\varphi_{0} = 0$,
        \item $w_t$ for $t = 0, \ldots, N-1$ is the $t$-th random noise vector,
		\item $\eta_t$ is an $m$-dimensional vector for $t = 0, \ldots, N-1$,
		\item $\varphi_t(w_0, \ldots, w_{t-1}) = \sum_{i=0}^{t-1}\varphi_{t, i}(w_i)$ for $t = 1, \ldots, N-1$ is an $m$-dimensional vector, and
		\item each function $\varphi_{t, i}$ belongs to the linear span of the basis elements $(\ee^\nu)_{\nu\in\mcal I}$, and thus has a representation as a linear combination $\varphi_{t, i}(\cdot) = \sum_{\nu\in\mcal I} \theta_{t, i}^\nu\ee^\nu(\cdot)$, $t = 1, \ldots, N-1$, $i = 0, \ldots, t-1$, where $\theta_{t, i}^\nu$ are matrices of coefficients of appropriate dimension.
	\end{itemize}
	Analogous to Fourier coefficients in harmonic analysis, we call the $\theta_{t, i}^\nu$ the $\nu$-th Fourier coefficient of the function $\varphi_{t, i}$. Therefore, whenever $|\mcal I| < \infty$ for every $t = 1, \ldots, N-1$, we have the finite representation
	\begin{equation}
	\label{e:thetadef}
		\varphi_t(w_0, \ldots, w_{t-1}) =
		\begin{bmatrix}
			\theta_{t, 0} & \theta_{t, 1} & \ldots & \theta_{t, t-1} & \pmb 0 & \ldots & \pmb 0
		\end{bmatrix}_{R^{m\times n|\mcal I|(N-1)}}
		\begin{bmatrix}
			\ee(w_{0})\\
			\ee(w_{1})\\
			\vdots\\
			\ee(w_{N-2})
		\end{bmatrix}_{\R^{n|\mcal I|(N-1)\times 1}},
	\end{equation}
   	where $\theta_{t, i}\in\R^{m\times n|\mcal I|}$, $\pmb 0\in \R^{m\times n|\mcal I|}$,
   	\begin{equation*}
   		\theta_{t,i}\Let\smat{\theta^1_{t,i}\,\; \cdots \,\, \theta^{|\mathcal I|}_{t,i} },\;\;\theta_{t, i}^\nu\in\R^{m\times n}, \quad {\rm and}\quad \ee(w_{i})\Let\smat{\ee^1(w_i)\\\vdots \\ \ee^{|\mathcal I|}(w_i)},\quad \fa i=0,1,\cdots, N-2.
   	\end{equation*}
	In this notation the policy~\eqref{e:genpolicy} can be written as
	\begin{equation}
	\label{e:Thetadef}
		 u  = \eta + \varphi( w ) = \eta +
		\begin{bmatrix}
			\pmb 0 & \pmb 0 & \cdots & \pmb 0\\
			\theta_{1, 0} & \pmb 0 & \cdots & \pmb 0\\
			\theta_{2, 0} & \theta_{2, 1} & \cdots & \pmb 0\\
			\vdots & \vdots & \ddots & \vdots\\
			\theta_{N-1, 0} & \theta_{N-1, 1} & \cdots & \theta_{N-1, N-2}
		\end{bmatrix}
		\begin{bmatrix}
			\ee(w_0)\\
			\ee(w_1)\\
			\vdots\\
			\ee(w_{N-2})
		\end{bmatrix}
		\teL \eta + \Theta \ee( w ),
	\end{equation}
	where $\Theta$ is now the matrix of Fourier coefficients having dimension $Nm\times\bigl(n(N-1)|\mcal I|\bigr)$. This Fourier coefficient matrix $\Theta$ and the vector $\eta$ play the role of the optimization parameters in our search for an optimal policy. Note that $\ee( w )$ does not include the noise vector $w_{N-1}$, and that $\Theta$ is strictly lower block triangular to enforce causality. In what follows, as a matter of notation, by $\Theta_t$ we shall denote the formal $t$-th block-row of the matrix $\Theta$ in~\eqref{e:Thetadef}, i.e., $\Theta_t \Let \begin{bmatrix} \theta_{t, 0} & \cdots & \theta_{t, t-1} & 0 & \cdots & 0\end{bmatrix}$, for $t=0, \cdots, N-1$, with $\Theta_0$ being the identically $0$ row. We make the following assumption:

	\begin{assumption}
	\label{a:sys}
		The sequence $(w_t)_{t\in\Nz}$ of noise vectors is i.i.d with $\Sigma = \EE\bigl[w_t w_t\transp\bigr]$.\AssmpEnd
	\end{assumption}

	So far we have not stipulated any boundedness properties on the elements of the vector space $\HH$. This means that the control policy elements may be unbounded maps. 
	First we stipulate the following structure on the control sets:

For a given $p\in[1, \infty]$, the control input vector $u_t$ is bounded in $p$-norm at each instant of time $t$, i.e., for $p\in[1, \infty]$ let $U_{\max}^{(p)} > 0$ be given, with
\begin{equation}
\begin{aligned}
\label{eqn:Utotal}
u_t\in\bar \U_p & \Let \bigl\{\xi\in\mathbb R^m\big| \norm{\xi}_p \le U_{\max}^{(p)}\bigr\} \quad \fa t\in\Nz,\quad\text{and}\\
\U_p & \Let \underbrace{\bar \U_p\times\ldots\times\bar  \U_p}_{N-\text{times}}.
\end{aligned}
\end{equation}
One could easily include more general constraint sets $\U_p$, for instance, to capture bounds on the rate of change of inputs.


	Our basic result is the next Theorem.

	\begin{theorem}
	\label{t:gen}
		Consider the system~\eqref{eq:system}. Suppose that Assumption~\ref{a:sys} holds, $\HH$ is finite-dimensional ($|\mcal I| < \infty$), and every component of the basis functions $\ee^\nu$ is bounded by $\mcal E > 0$ in absolute value. Then the problem~\eqref{eq:problem1} under the policy~\eqref{e:genpolicy} and control sets \eqref{eqn:Utotal} for $p\in[1, \infty]$ is convex with respect to the decision variables $(\eta, \Theta)$ defined in~\eqref{e:Thetadef}. For $p = 1, 2$, and $\infty$ it admits convex tractable versions with tighter domains of $(\eta, \Theta)$, given by
		\begin{equation}
		\label{e:probspan}
		\begin{aligned}
			\underset{(\eta, \Theta)}{\text{minimize}} \quad & \trace\Bigl(\Theta\transp\bigl({ B }\transp Q  B + R \bigr)\Theta \Sigma_{\ee}\Bigr) + 2\trace\Bigl(\Theta\transp{ B }\transp Q  D \Sigma_{\ee}'\Bigr) + \eta\transp\bigl({ B }\transp Q  B  +  R \bigr)\eta\\
			& \quad + 2\bigl(\xz\transp{ A }\transp Q  B \eta + \eta\transp{ B }\transp Q  D \mu_w + x_0\transp A \transp Q  B \Theta\mu_\ee\bigr)\\
			& \quad + 2\eta\transp\bigl( R  +  B \transp Q  B \bigr)\Theta\mu_\ee + c\\
			\text{subject to }\quad & \text{$\Theta$ strictly lower block triangular as in~\eqref{e:Thetadef}},\\
			& \begin{cases}
			p=1: &	\norm{\eta_t}_1 + \mcal Et \norm{\Theta_t}_1 \le U_{\max}^{(1)},\quad \fa t = 0, 1, \ldots, N-1,\\
			p=\infty: &	\norm{\eta_t}_\infty + \mcal E\norm{\Theta_t}_\infty \le U_{\max}^{(\infty)}, ,\quad \fa t = 0, 1, \ldots, N-1, \\
			p=2: &	\norm{\begin{bmatrix}\eta_t & \Theta_t\end{bmatrix}}_2\sqrt{1+\mcal Et} \le U_{\max}^{(2)}, \quad \fa t = 0, 1, \ldots, N-1,
			\end{cases}
		\end{aligned}
		\end{equation}
		where
		\begin{align*}
			\Sigma_{\ee} & \Let \EE\bigl[\ee( w )\ee( w )\transp\bigr], &  \Sigma_{\ee}' & \Let \EE\bigl[ w \ee( w )\transp\bigr],\\
			\mu_w & \Let \EE[ w ], \quad \quad \mu_\ee \Let\EE[\ee( w )], & c & \Let x_0\transp{ A }\transp Q  A x_0 + 2\xz\transp{ A }\transp Q  D  \mu + \trace\bigl({ D }\transp Q  D \Sigma_{ w }\bigr).
		\end{align*}
	\end{theorem}


	\begin{proof}[Proof of Theorem~\ref{t:gen}]
		It is easy to see that  $x\transp Q x + u\transp R u$ is convex nondecreasing, and both $x$ and $u$ are affine functions of the design parameters $(\eta,\Theta)$ for any realization of the noise $w$. Hence, $V_0$ is convex in $(\eta, \Theta)$ since taking expectations of a convex function retains convexity \cite[Section 3.2]{ref:boyd04}. Moreover, the control constraint sets in \eqref{eqn:Utotal} are convex in $(\eta, \Theta)$. This settles the first claim.

		The objective function~\eqref{e:objfn} is given by
		\begin{align*}
			\EE_{\xz} & \bigl[\bigl( A  x_0 +  B  u  +  D  w \bigr)\transp Q \bigl( A  x_0 +  B  u  +  D  w \bigr)\bigr] + \EE_{\xz}\bigl[{ u }\transp R  u \bigr]\\
			& = \EE_{\xz}\bigl[\bigl( A  x_0 +  B (\eta + \Theta \ee ( w )) +  D  w \bigr)\transp  Q \bigl( A  x_0 +  B (\eta + \Theta \ee ( w )) +  D  w \bigr)\bigr]\\
			& \quad + \EE_{\xz}\bigl[(\eta + \Theta \ee ( w ))\transp R (\eta + \Theta \ee ( w ))\bigr]\\
			& = x_0\transp{ A }\transp Q  A x_0 + 2\xz\transp{ A }\transp Q  B \eta + \eta\transp\bigl({ B }\transp Q  B  +  R \bigr)\eta \\
			& \quad + 2\bigl( A \xz +  B \eta\bigr)\transp Q \bigl( B \Theta\EE_{\xz}[ \ee ( w )] +  D \EE_{\xz}[ w ]\bigr)\\
			& \quad + \EE_{\xz}\bigl[\bigl( B \Theta \ee ( w ) +  D  w \bigr)\transp Q \bigl( B \Theta \ee ( w ) +  D  w \bigr)\bigr] + \EE_{\xz}\bigl[(\Theta \ee ( w ))\transp R \Theta \ee ( w )\bigr]\\
			& = x_0\transp{ A }\transp Q  A x_0 + 2\xz\transp{ A }\transp Q  B \eta + \eta\transp\bigl({ B }\transp Q  B  +  R \bigr)\eta + 2\eta\transp R \Theta\EE[ \ee ( w )]\\
			& \quad +  2\bigl( A \xz +  B \eta\bigr)\transp Q \bigl(  D \EE_{\xz}[ w ] +  B \Theta\EE[ \ee ( w )]\bigr) + \trace\Bigl({ D }\transp Q  D \EE_{\xz}\bigl[ w \transp\bigr]\Bigr)\\
			& \quad + \trace\Bigl(\Theta\transp\bigl({ B }\transp Q  B + R \bigr)\Theta\EE_{\xz}\bigl[ \ee ( w ) \ee ( w )\transp\bigr]\Bigr) + 2\trace\Bigl(\Theta\transp{ B }\transp Q  D \EE_{\xz}\bigl[ w  \ee ( w )\transp\bigr]\Bigr).
		\end{align*}
		Incorporating the definitions $\Sigma_{\ee}$, $\Sigma_{\ee}'$, $\mu_w$, $\mu_\ee$, and $c$, the right-hand side above equals
		\begin{align*}
			& \trace\Bigl(\Theta\transp\bigl({ B }\transp Q  B + R \bigr)\Theta\Sigma_{\ee}\Bigr) + 2\trace\Bigl(\Theta\transp{ B }\transp Q  D \Sigma_{\ee}'\Bigr) + \eta\transp\bigl({ B }\transp Q  B  +  R \bigr)\eta\\
			& \quad + 2\bigl(\xz\transp{ A }\transp Q  B \eta + \eta\transp{ B }\transp Q  D \mu_w + x_0\transp A \transp Q  B \Theta\mu_\ee\bigr) + 2\eta\transp\bigl( R  +  B \transp Q  B \bigr)\Theta\mu_\ee \\
			& \quad + \bigl(x_0\transp{ A }\transp Q  A x_0 + 2\xz\transp{ A }\transp Q  D  \mu_w + \trace\bigl({ D }\transp Q  D \Sigma_{ w }\bigr)\bigr)\\
			& = \trace\Bigl(\Theta\transp\bigl({ B }\transp Q  B + R \bigr)\Theta\Sigma_{\ee}\Bigr) + 2\trace\Bigl(\Theta\transp{ B }\transp Q  D \Sigma_{\ee}'\Bigr) + \eta\transp\bigl({ B }\transp Q  B  +  R \bigr)\eta\\
			& \quad + 2\bigl(\xz\transp{ A }\transp Q  B \eta + \eta\transp{ B }\transp Q  D \mu_w + x_0\transp A \transp Q  B \Theta\mu_\ee\bigr) + 2\eta\transp\bigl( R  +  B \transp Q  B \bigr)\Theta\mu_\ee + c.
		\end{align*}
		Since the matrix $\Sigma_{\ee}$ is positive semidefinite, it can be expressed as a finite nonnegative linear combination of matrices of the type $\sigma\sigma\transp$, for vectors $\sigma$ of appropriate dimension~\cite[Theorem~1.10]{ref:bermanCPM}. Accordingly, if $\Sigma_{\ee} = \sum_{i=1}^k \sigma_i\sigma_i\transp$, then
		\begin{align*}
			\trace\Bigl(\Theta\transp\bigl({ B }\transp Q  B + R \bigr)\Theta\Sigma_{\ee}\Bigr) &= \sum_{i=1}^k\trace\bigl(\Theta\transp\bigl({ B }\transp Q  B + R \bigr)\Theta\sigma_i\sigma_i\transp\bigr)\\
			& = \sum_{i=1}^k \Bigl(\sigma_i\transp\Theta\transp\bigl({ B }\transp Q  B + R \bigr)\Theta\sigma_i\Bigr).
		\end{align*}
		Defining $\wh\Theta_i \Let \Theta \sigma_i$ and adjoining these equalities to the constraints of the optimization program~\eqref{e:probspan}, we arrive at the optimization program
		\begin{equation}
		\label{e:genqp}
		\begin{aligned}
			\underset{(\Theta, \wh\Theta_1, \ldots, \wh\Theta_k)}{\text{minimize}} \quad & \sum_{i=1}^k\wh\Theta_i\transp\bigl({ B }\transp  Q  B  +  R \bigr)\wh\Theta_i + 2\trace\Bigl(\Theta\transp{ B }\transp Q  D \Sigma_{\ee}'\Bigr) + \eta\transp\bigl({ B }\transp Q  B  +  R \bigr)\eta \\
			& \qquad + 2\bigl(\xz\transp{ A }\transp Q  B \eta + \eta\transp{ B }\transp Q  D \mu + x_0\transp A \transp Q  B \Theta\mu^\ee\bigr)\\
			& \qquad + 2\eta\transp\bigl( R  +  B \transp Q  B \bigr)\Theta\mu^\ee + c\\
			\text{subject to}\quad & \Theta \text{ strictly lower block triangular as in~\eqref{e:Thetadef}},\\
			& \wh\Theta_i = \Theta\sigma_i \quad\text{for all }i = 1, \ldots, k.
		\end{aligned}
		\end{equation}
		We see immediately that~\eqref{e:genqp} is a convex program in the parameters $\eta$, $\Theta$ and $\wh\Theta_i$, and is equivalent to the cost in \eqref{e:probspan}.

		It only remains to consider the last constraint in \eqref{e:probspan}. First we consider the cases of $p = 1, \infty$. Using the notation above, an application of the triangle inequality immediately shows that the constraints can be written as
		\begin{equation}
		\label{e:genqpconstr}
		\begin{aligned}
    	\begin{cases}
    	p = 1: & \norm{\eta_t}_1 + \mcal E t\norm{\Theta_t}_1 \le U_{\max}^{(1)},\quad \fa t = 0, 1,\ldots, N-1,\\
    	p=\infty: &		\norm{\eta_t}_\infty + \mcal E\norm{\Theta_t}_\infty \le U_{\max}^{(\infty)},\quad \fa t = 0, 1,\ldots, N-1.
    	\end{cases}
		\end{aligned}
		\end{equation}
		It follows that the objective function in \eqref{e:genqp} is quadratic and the constraints in \eqref{e:genqp}-\eqref{e:genqpconstr} are affine in the optimization parameters $\eta$, $\Theta$, and $\wh\Theta$. As such, for $p = 1, \infty$ our problem is a quadratic program.

		For the case of $p = 2$, note that $\eta_t + \Theta_t \ee( w ) = \begin{bmatrix}\eta_t & \Theta_t\end{bmatrix}\begin{bmatrix}1\\ \ee( w )\end{bmatrix}$, and by definition of $\mcal E$ it is clear that $\norm{\begin{bmatrix}\eta_t & \Theta_t\end{bmatrix}\begin{bmatrix}1\\ \ee( w )\end{bmatrix}}_2 \le \norm{\begin{bmatrix}\eta_t & \Theta_t\end{bmatrix}}_2 \sqrt{1 + \mcal Et}$. This immediately translates to $\norm{\begin{bmatrix}\eta_t & \Theta_t\end{bmatrix}}_2 \sqrt{1+\mcal Et}\le U_{\max}^{(2)}$, which is the third constraint in Problem \ref{e:probspan} and it is a quadratic constraint in the optimization parameters $(\eta, \Theta)$. Therefore, for $p=2$ our problem is a quadratically constrained quadratic program.
	\end{proof}

	The optimization problem \eqref{e:probspan} simplifies if we assume that $\mu^\ee = \EE[ \ee ( w )] = 0$. Note that $\EE[ \ee ( w )] = 0$ if and only if $\EE\bigl[\ee_{t, i}^\nu(w_{t, i})\bigr] = 0$ for all $\nu\in\mcal I$. At an intuitive level this translates to the condition that the functions $\ee_{t,i}^\nu\in\HH$ should be ``centered'' with respect to the random variables $w_{t, i}$. In particular, this simply means that for noise distributions that are symmetric about $0$, the functions $\ee^\nu$ should be centered at $0$ and be antisymmetric. For example, if the noise is Gaussian with mean $0$ and diagonal covariance matrix (uncorrelated components), each component of the functions $\ee^\nu$ should be an odd function.

	The matrices $\Sigma_{\ee}$, $\Sigma_{\ee}'$, the vector $v$, and the number $c$ in Theorem~\ref{t:gen} are all constants independent of $x_0$, and can be computed off-line. As such, even if closed-form expressions for the entries of the matrices do not exist, they can be numerically computed to desired precision. The optimization problem~\eqref{e:probspan} is a quadratic program~\cite[p.~152]{ref:boyd04} for $p = 1, \infty$, and a quadratically constrained quadratic program~\cite[p.~152]{ref:boyd04} for $p = 2$, in the optimization parameters $\bigl\{\eta, \Theta, \bigl\{\wh\Theta_i, i=1, \ldots, k\bigr\}\bigr\}$, and can be easily coded in standard software packages such as \texttt{cvx}~\cite{ref:boydCVX} or \texttt{YALMIP} \cite{YALMIP}. Note that the optimization problem~\eqref{e:probspan} is always feasible (simply set $\Theta = 0$ and $\eta = 0$ to see this). This is not a surprise, since there are no constraints on the state, and by construction $0\in\mathbb U$. Finally, note that the third constraint in Problem \eqref{e:probspan} for various values of $p$, is a result of robustly satisfying the constraints posed by the various control sets \eqref{eqn:Utotal} for any realization of the noise $w$.

	In general, the total number of decision variables in the optimization program \eqref{e:probspan} is $mN\bigl(1+\tfrac{1}{2}n(N-1)|\mcal I|\bigr)$. The number of decision variables can be substantially reduced, e.g., by choosing $\HH$ to be $1$-dimensional, or by fixing certain (block) elements of the Fourier coefficient matrix $\Theta$ to $0$.


\section{Various Cases of Constrained Controls}
\label{s:if}

We examine in this section several special cases of Theorem~\ref{t:gen} under various restrictions on the classes of noise and control inputs.


	\subsection{Bounded controls, unbounded noise, and $p=\infty$}\label{sec:bdduubddw}
		Let the noise take values in $\R^n$. We provide tractable convex programs to design a policy that by construction respects the control constraint sets \eqref{eqn:Utotal}, with $p=\infty$. Starting from~\eqref{e:genpolicy} let
		\begin{equation}
		\label{eq:auginputbdd}
			u = \eta + \Theta\varphi(w),
		\end{equation}
		where
		\begin{itemize}[label=\textbullet, leftmargin=*]
			\item $\varphi(w)\Let \smat{\varphi_0\\ \varphi_1(w_0)\\ \vdots \\ \varphi_{N-1}(w_0, \ldots, {w}_{N-2})}$,
			\item $\varphi_0 = 0$, $\varphi_t(w_0, \ldots, w_{t-1}) = \sum_{j=0}^{t-1} \theta_t^j\varphi_{t, j}(w_j)$ for $t = 1, \ldots, N-2$, and
			\item $\varphi_{t, j}(w_j) = \bigl[\wt\varphi(w_{j, 1}), \ldots, \wt\varphi(w_{j,n})\bigr]\transp$ for some function $\wt\varphi$ such that $\sup\limits_{s\in\R}\wt\varphi(s) = \phi_{\max} < \infty$, and $\varphi_{t, j}:\W\to\U_\infty$.
		\end{itemize}
In other words, we saturate the measurements that we obtain from the noise input vector before inserting them into our control vector. This way we allow that the noise distribution is supported over the entire $\R^n$, which is an advantage over other approaches \cite{BertsimasBrown-07,ref:goulart06}. Moreover, the choice of the component saturation function $\wt\varphi$ is left open as long as the noise sequence satisfies Assumption~\ref{a:sys}. For example, we can accommodate standard saturation, piecewise linear, and sigmoidal functions to name a few.

			Our choice of saturating the measurement from the noise vectors, as we shall see below, renders the resulting optimization problem tractable as opposed to calculating the entire control input vector $  u$ and then saturating it a posteriori; one can see that the latter approach tends to lead to an intractable optimization problem. Note also that the choice of control inputs in \eqref{eq:auginputbdd} yields a possibly non-Markovian feedback. 

		\begin{corollary}
		\label{c:sat}
			Consider the system~\eqref{eq:system}. Suppose that Assumption~\ref{a:sys} holds, and $\EE[ \ee(  w)] = 0$ with $ \ee(  w) = \varphi( w)$, where $\varphi$ is defined in~\eqref{eq:auginputbdd}. Then for $p = \infty$ the problem~\eqref{eq:problem1} under the control policy~\eqref{eq:auginputbdd} is a convex optimization program with respect to the decision variables $(\eta, \Theta)$, given by
			\begin{equation}
			\label{eq:problem2}
			\begin{aligned}
				\underset{(\eta, \Theta)}{\text{minimize}} \quad& \tr{\Theta\transp \bigl( R +  B \transp  Q  B \bigr)\Theta\Gamma_1} +2\tr{  D  Q  B \Theta\Gamma_2}\\
				& + \eta\transp \bigl( R +  B \transp  Q  B \bigr)\eta +b\transp \eta+c \\
				\text{subject to}\quad& \max\limits_{i=1,\cdots, m} \left(|\eta_{t,i}| + \norm{\Theta_{t,i}}_1\phi_{\rm max}\right) \le U_{\max}^{(\infty)},\quad t=0, \ldots, N-1,\\
				& \text{and $\Theta$ strictly lower block triangular as in~\eqref{e:Thetadef}},
			\end{aligned}
			\end{equation}
			where $\eta_{t,i}$ and $\Theta_{t,i}$ are the $i$-th rows of $\eta_t$ and $\Theta_t$, respectively,
			\begin{align*}
				c &= x_0\transp  A  Q  A  x_0 + \tr{ D \transp  Q  D  \Sigma_{\bar w}}, \\
				b &= 2 B \transp  Q  A  x_0, \\
				\Gamma_1 &= \mathrm{diag}\bigl\{\EE\bigl[\varphi_0(w_0)\varphi_0(w_0)\transp\bigr],\cdots,\EE\bigl[\varphi_{N-1}(w_{N-1})\varphi_{N-1}(w_{N-1})\transp\bigr]\bigr\},\\
				\Gamma_2 &= \mathrm{diag}\bigl\{\EE\bigl[\varphi_0(w_0) w_0\transp\bigr],\cdots,\EE\bigl[\varphi_{N-1}(w_{N-1})w_{N-1}\transp\bigr]\bigr\}.
			\end{align*}
			The resulting policy is guaranteed to satisfy the control constraint set \eqref{eqn:Utotal} for $p = \infty$.
		\end{corollary}

		A complete proof may be found in~\cite{ref:petercdc09}; it proceeds along the lines of the proof of Theorem~\ref{t:gen}. Note that the program~\eqref{eq:problem2} exactly solves~\eqref{eq:problem1} under the policy~\eqref{eq:auginputbdd} and is neither a restriction nor a relaxation.

			Problem~\eqref{eq:problem2} is a quadratic program in the optimization parameters $(\eta, \Theta)$ (see the discussion following Theorem \ref{t:gen}). The matrices $\Gamma_1$ and $\Gamma_2$ capture the statistics of the noise in the presence of the functions $\varphi$ and can be computed numerically \emph{off-line} using Monte Carlo techniques \cite[Section 3.2]{RobertCasella-04}. This method will be utilized in the examples in Section \ref{sec:examples}. However, in some instances it is actually possible to compute these matrices in closed form; this is shown in the next three examples.

		\begin{example}
		\label{ex:sigmoids}
			Let us consider~\eqref{eq:system} when the noise process $(w_t)_{t\in\Nz}$ is an i.i.d sequence of Gaussian random vectors of mean $0$ and covariance $\Sigma$ and standard sigmoidal policy functions $\wt\varphi$, i.e., $\wt\varphi(t) \Let t/\sqrt{1+t^2}$. Assume further that the components of $w_t$ are mutually independent, which implies that $\Sigma$ is a diagonal matrix $\diag\{\sigma_1^2, \ldots, \sigma_n^2\}$.  Then from the identities in Fact \ref{facts} in \secref{s:facts}, we have for $i=1, \ldots, n$ and $j=0, \ldots, N-1$,
			\begin{align*}
				\EE\bigl[\wt\varphi(w_{j, i})^2\bigr] & = \frac{1}{\sqrt{2\pi}\sigma_i}\int_{-\infty}^\infty \wt\varphi(t)^2 \epower{-\frac{t^2}{2\sigma_i^2}} \mrm dt = 2\cdot\frac{1}{\sqrt{2\pi}\sigma_i}\int_{0}^\infty \frac{t^2}{1+t^2} \epower{-\frac{t^2}{2\sigma_i^2}}\\
				& = \sqrt{2\pi}\sigma_i - \pi\epower{-\frac{1}{2\sigma_i^2}}\erfc\Bigl(\frac{1}{\sqrt 2\sigma_i}\Bigr).
			\end{align*}
			This shows that the matrix $\Gamma_1$ in Corollary~\ref{c:sat} is $\diag\{\Sigma', \ldots, \Sigma'\}$, where
			\[
				\Sigma' \Let \diag\left\{\sqrt{2\pi}\sigma_1 - \pi\epower{-\frac{1}{2\sigma_1^2}}\erfc\Bigl(\frac{1}{\sqrt 2\sigma_1}\Bigr), \ldots, \sqrt{2\pi}\sigma_n - \pi\epower{-\frac{1}{2\sigma_n^2}}\erfc\Bigl(\frac{1}{\sqrt 2\sigma_n}\Bigr)\right\}.
			\]
			Similarly, since
			\begin{align*}
				\EE\bigl[\wt\varphi(w_{j,i}) w_{j,i}\bigr] & = \frac{1}{\sqrt{2\pi}\sigma_i}\int_{-\infty}^\infty t\wt\varphi(t) \epower{-\frac{t^2}{2\sigma_i^2}} \mrm dt = 2\cdot\frac{1}{\sqrt{2\pi}\sigma_i}\int_{-\infty}^\infty \frac{t^2}{\sqrt{1+t^2}} \epower{-\frac{t^2}{2\sigma_i}} \mrm dt\\
				& = \frac{\sigma_i}{\sqrt 2}U\Bigl(\frac{1}{2}, 0, \frac{1}{2\sigma_i^2}\Bigr),
			\end{align*}
			where $U$ is the confluent hypergeometric function (defined in the Appendix), the matrix $\Gamma_2$ in Corollary~\ref{c:sat} is $\diag\{\Sigma'', \ldots, \Sigma''\}$, where
			\[
				\Sigma'' \Let \diag\left\{\frac{\sigma_1}{\sqrt 2}U\Bigl(\frac{1}{2}, 0, \frac{1}{2\sigma_1^2}\Bigr), \ldots, \frac{\sigma_n}{\sqrt 2}U\Bigl(\frac{1}{2}, 0, \frac{1}{2\sigma_n^2}\Bigr)\right\}.
			\]
			Therefore, given the system~\eqref{eq:system}, the control policy~\eqref{eq:bentalpolicy}, and the description of the noise input as above, the matrices $\Gamma_1$ and $\Gamma_2$ derived above complete the set of hypotheses of Corollary~\ref{c:sat}. The problem~\eqref{eq:problem} can now be solved as the quadratic program~\eqref{eq:problem2}.\ExEnd
		\end{example}

		\begin{example}
		\label{ex:sigmoidsgen}
			Consider the setting of Example~\ref{ex:sigmoids} (with $\wt\varphi$ a standard sigmoid) under the assumption that $\Sigma$ is a not necessarily diagonal matrix. To wit, the components of $w_t$ may be correlated at each time $t\in\Nz$; however, the random vector sequence $(w_t)_{t\in\Nz}$ is assumed to be i.i.d. This is equivalent to the knowledge of the correlations between the random variables $\bigl\{w_{t, i}\big|i=1, \ldots, n\bigr\}$, which are constant over $t$. Then $\EE[\varphi(\bar w)\varphi(\bar w)\transp]$ is a block diagonal matrix. Indeed, we have with $\Sigma_{i, j} \Let \begin{bmatrix}\sigma_i^2 & \rho_{i, j}^2\\\rho_{i, j}^2 & \sigma_j^2\end{bmatrix}$,
			\begin{align*}
				\EE\bigl[ & \wt\varphi(w_{t, i})\wt\varphi(w_{t, j})\bigr] \\
				& = \frac{1}{\sqrt{2\pi\det{\Sigma_{i, j}}}}\iint_{\R^2} \frac{t_1 t_2}{\sqrt{(1+t_1^2)(1+t_2^2)}} \exp\biggl(-\frac{1}{2}\begin{bmatrix}t_1 & t_2\end{bmatrix}\Sigma_{i, j}^{-1}\begin{bmatrix}t_1\\t_2\end{bmatrix}\biggr) \;\mrm dt_1\mrm dt_2,
			\end{align*}
			and
			\begin{align*}
				\EE\bigl[ & \wt\varphi(w_{t,i}) w_{t,j}\bigr] & = \frac{1}{\sqrt{2\pi\det\Sigma_{i, j}}}\iint_{\R^2} \frac{t_1t_2}{\sqrt{1+t_1^2}} \exp\biggl(-\frac{1}{2}\begin{bmatrix} t_1 & t_2\end{bmatrix}\Sigma_{i, j}^{-1} \begin{bmatrix}t_1\\t_2\end{bmatrix}\biggr) \;\mrm dt_1\mrm dt_2.
			\end{align*}
			Note that the computations of the integrals above can be carried out off-line. We define the matrices $\Sigma_t$ and $\Sigma_t'$ with the $(i,j)$-th entry of $\Sigma_t$ being $\EE\bigl[\wt\varphi(w_{t,i})\wt\varphi(w_{t, j})\bigr]$ and the $(i, j)$-th entry of $\Sigma_t'$ being $\EE\bigl[\wt\varphi(w_{t, i})w_{t, j}\bigr]$, and it follows that the matrices $\Gamma_1 = \diag\bigl\{\Sigma_0, \ldots, \Sigma_{N-2}\bigr\}$, and $\Gamma_2 = \diag\bigl\{\Sigma_0', \ldots, \Sigma_{N-2}'\bigr\}$. \ExEnd
		\end{example}

		\begin{example}
		\label{ex:sat}
			Consider the system~\eqref{eq:system} as in Example~\ref{ex:sigmoids}, and with $\wt\varphi$ the standard saturation function defined as $\wt\varphi(t) = \sgn(t)\min\{|t|, 1\}$. From Corollary~\ref{c:sat} we have for $i=1, \ldots, n$ and $j=0, \ldots, N-1$, using the identities in Fact \ref{facts} in \secref{s:facts},
			\begin{align*}
				\EE\bigl[\wt\varphi(w_{j,i})^2\bigr] & = \frac{1}{\sqrt{2\pi}\sigma_i}\int_{-\infty}^\infty \wt\varphi(t)^2 \epower{-\frac{t^2}{2\sigma_i^2}} \mrm dt\\
				& = \frac{2}{\sqrt{2\pi}\sigma_i}\int_0^1 t^2\epower{-\frac{t^2}{2\sigma_i^2}} \mrm dt + \frac{2}{\sqrt{2\pi}\sigma_i}\int_1^\infty \epower{-\frac{t^2}{2\sigma_i^2}} \mrm dt\\
				& = \sqrt{2\pi}\sigma_i^3\erf\Bigl(\frac{1}{\sqrt 2\sigma_i}\Bigr) - 2\sigma_i^2\epower{-\frac{1}{2\sigma_i^2}} + 1 + \erf\Bigl(\frac{1}{\sqrt 2\sigma_i}\Bigr)\\
				& \teL \xi_i' \text{ (say)},
			\end{align*}
			and
			\begin{align*}
				\EE\bigl[\wt\varphi(w_{j,i}) w_{j,i}\bigr] & = \frac{1}{\sqrt{2\pi}\sigma_i}\int_{-\infty}^\infty t \wt\varphi(t) \epower{-\frac{t^2}{2\sigma_i^2}} \mrm dt\\
				& = \frac{2}{\sqrt{2\pi}\sigma_i}\int_0^1 t^2\epower{-\frac{t^2}{2\sigma_i^2}} \mrm dt + \frac{2}{\sqrt{2\pi}\sigma_i}\int_1^\infty t\epower{-\frac{t^2}{2\sigma_i^2}} \mrm dt\\
				& = \sqrt{2\pi}\sigma_i^3\erf\Bigl(\frac{1}{\sqrt 2\sigma_i}\Bigr) - 2\sigma_i^2\epower{-\frac{1}{2\sigma_i^2}} + \sqrt{\frac{2}{\pi}}\sigma_i\Gammaf(2\sigma_i^2, 1)\\
				& \teL \xi_i'' \text{ (say)}.
			\end{align*}
			Therefore, in this case the matrix $\Gamma_1$ in Corollary~\ref{c:sat} is $\diag\{\Sigma', \ldots, \Sigma'\}$ with $\Sigma' \Let \diag\{\xi_1', \ldots, \xi_n'\}$, and the matrix $\Gamma_2$ is $\diag\{\Sigma'', \ldots, \Sigma''\}$ with $\Sigma'' \Let \diag\{\xi_1'', \ldots, \xi_n''\}$. These information complete the set of hypotheses of Corollary~\ref{c:sat}, and the problem~\eqref{eq:problem} can now be solved as a quadratic program~\eqref{eq:problem2}.\ExEnd
		\end{example}


	\subsection{Bounded controls, bounded noise, and $p=2$}\label{sec:bcbn}
		In this subsection we specialize to the case of the noise being drawn from a compact subset of $\R^n$, and the control inputs set $\U_2$. We make the following assumption:
		\begin{assumption}
		\label{a:w}
			The noise takes values in a compact set $\W \subset \R^n$.\AssmpEnd
		\end{assumption}
		Under Assumption~\ref{a:w} Hilbert space techniques may be effectively employed in our basic controller synthesis framework of Section \ref{s:mainres} in the following way.
		Let $(\HH, \inprod{\cdot}{\cdot}_{\HH})$ be a separable Hilbert space of measurable maps $\ee:\W\to\U_2$ supported on the compact set $\W$. The inner product is defined as $\inprod{\varphi_1}{\varphi_2}_{\HH} \Let \sum_{i=1}^n \inprod{\varphi_{1, i}}{\varphi_{2, i}}$ where $\inprod{\cdot}{\cdot}$ is the standard inner product on real-valued functions on $\W$. Fix a complete orthonormal basis $(\ee^\nu)_{\nu\in\mcal I}\subset\HH$. Since $\HH$ is separable, the set $\mcal I$ is at most countable. Just as in~\eqref{e:genpolicy} we let our candidate control policies be of the form
		\begin{equation}
		\label{e:hspolicy}
			  u =
			\begin{bmatrix}
				\eta_0\\
				\eta_1\\
				\vdots\\
				\eta_{N-1}
			\end{bmatrix} +
			\begin{bmatrix}
				\mathbf 0 & \mathbf 0 & \cdots & \mathbf 0\\
				\theta_{1, 0} & \mathbf 0 & \cdots & \mathbf 0\\
				\theta_{2, 0} & \theta_{2, 1} & \cdots & \mathbf 0\\
				\vdots & \vdots & \ddots & \vdots\\
				\theta_{N-1, 0} & \theta_{N-1, 1} & \cdots & \theta_{N-1, N-2}
			\end{bmatrix}
			\begin{bmatrix}
				 \ee(w_0)\\
				 \ee(w_1)\\
				\vdots\\
				 \ee(w_{N-2})
			\end{bmatrix}
			\teL \eta + \Theta \ee( w),
		\end{equation}
		where the vector $ \ee(\cdot)$ is the formal vector formed by concatenating the (ordered) basis elements $(\ee^\nu)_{\nu\in\mcal I}$, the various $\theta$-s are formal matrices as in Section \ref{s:mainres}, and $\eta_t$ is an $m$-dimensional vector for $t = 0, \ldots, N-1$. This takes us back to the setting of Section \ref{s:mainres}.

		The following Corollary illustrates the technique explained above; its proof will only be sketched---it is similar to the proof of Theorem \ref{t:gen}. Note that for finite-dimensional Hilbert spaces, depending on the choice of the orthonormal basis, the matrix $\Theta$ may have complex or real entries.

		\begin{corollary}
		\label{c:hspolicy}
			Consider the system~\eqref{eq:system}. Suppose that Assumptions \ref{a:sys} and \ref{a:w} hold. Then for $p = 2$ and corresponding control set $\U_2$ problem~\eqref{eq:problem1} under the policy~\eqref{e:hspolicy} admits convex tractable reformulation with tighter domains of the decision variables $(\eta, \Theta)$ defined in~\eqref{e:hspolicy}, and is equivalent to the following program:
			\begin{equation}
			\label{e:probhs}
			\begin{aligned}
				& \text{the minimization problem~\eqref{e:probspan}}\\
				& \text{subject to}\quad \norm{\eta_t} + \sqrt{N-1}\norm{\Theta_t} \le U_{\max}^{(2)}, \quad\text{for }t=0, \ldots, N-1,\\
				& \qquad\qquad\qquad\text{and $\Theta$ strictly lower block triangular as in~\eqref{e:Thetadef}}.
			\end{aligned}
			\end{equation}
			Moreover, if $\hat\HH$ is a finite-dimensional subspace of $\HH$ spanned by $(\ee^\nu)_{\nu\in\mcal J}$ for some finite $\mcal J\subset\mcal I$, then the problem~\eqref{e:probhs} admits a reformulation as a quadratically constrained quadratic program with respect to the new decision variables $\bigl(\eta, \hat\Theta\bigr)$ corresponding to $\hat\HH$, given by
			\begin{equation}
			\label{e:probhsrel}
			\begin{aligned}
				& \text{the minimization problem~\eqref{e:probspan}}\\
				&\text{subject to}\quad \norm{\begin{bmatrix}\eta_t & \hat\Theta_t\end{bmatrix}} \le U_{\max}^{(2)}/\sqrt N \quad\text{for }t=0, \ldots, N-1,\\
				& \qquad\qquad\qquad\text{and $\Theta$ strictly lower block triangular as in~\eqref{e:Thetadef}},
			\end{aligned}
			\end{equation}
			where the vector ${\hat\ee}(\cdot)$ is the vector formed by concatenating the (ordered) basis elements $(\ee^\nu)_{\nu\in\mcal J}$, $ {\hat\ee}(  w) \Let \bigl[ {\hat\ee}(w_0)\transp, \ldots,  {\hat\ee}(w_{N-2})\transp\bigr]\transp$, $\hat\Sigma_{\ee} \Let \EE\bigl[ {\hat\ee}(  w) {\hat\ee}(  w)\transp\bigr]$, $\hat\Sigma_{\ee}' \Let \EE\bigl[  w {\hat\ee}(  w)\transp\bigr]$. In both the above cases the resulting policies are guaranteed to satisfy the control constraint set \eqref{eqn:Utotal} for $p = 2$.
		\end{corollary}
		\begin{proof}
			(Sketch.) Evaluating the objective function in~\eqref{eq:problem1} gives the objective function in~\eqref{e:probspan}. Recall that $\Theta_t$ is the $t$-th block row of the formal matrix $\Theta$, and $\Theta_{t, i}$ is the $i$th sub-row of the block row $\Theta_t$, where $t=0, \ldots, N-1$ and $i=1, \ldots, n$. Applying the triangle inequality for any $t = 0, \ldots, N-1$, we get
			\begin{align*}
				\norm{\eta_t + \Theta_t \ee(  w)} & \le \norm{\eta_t} + \norm{\Theta_t \ee(  w)} = \norm{\eta_t} + \sqrt{\inprod{\Theta_t \ee(  w)}{\Theta_t \ee(  w)}_{\HH}}\\
				& = \norm{\eta_t} + \sqrt{\sum_{i=1}^n \inprod{\Theta_{t, i} \ee(  w)}{\Theta_{t, i} \ee(  w)}} = \norm{\eta_t} + \sqrt{(N-1)\sum_{i=1}^n \norm{\Theta_{t, i}}^2}\\
				& = \norm{\eta_t} + \sqrt{N-1}\norm{\Theta_t}
			\end{align*}
			by orthogonality of the basis elements $(\ee^\nu)_{\nu\in\mcal I}$. The right-hand side of the last equality appears as the constraint in~\eqref{e:probhs}.

For the finite-dimensional case~\eqref{e:probhsrel}, we note that the objective function is identical to the one in~\eqref{e:probhs}, and the constraint in~\eqref{e:probhsrel} follows from the fact that $\norm{\eta_t + \hat\Theta_t {\hat\ee}(  w)} = \norm{\begin{bmatrix}\eta_t & \hat\Theta_t\end{bmatrix}\begin{bmatrix}1\\ {\hat\ee}(  w)\end{bmatrix}}$, and $\norm{\begin{bmatrix}1\\ {\hat\ee}(  w)\end{bmatrix}} = \sqrt{1+\sum_{i=0}^{N-2}\inprod{ {\hat\ee}(w_i)}{ {\hat\ee}(w_i)}} = \sqrt{N}$. This leads to a quadratically constrained quadratic program in the finite- dimensional decision variables $\bigl(\eta, \hat\Theta\bigr)$.
		\end{proof}Let us illustrate the usage of Corollary \ref{c:hspolicy} through the following example.
		\begin{example}
		\label{ex:unif}
			Consider the system~\eqref{eq:system}, and suppose that the $n$ components of the noise vector $w_t$ are independent uniform random variables taking values in $[-a, a]$ for some $a > 1$. Therefore, $\W = [-a, a]^n$. 
			It is a standard fact in Fourier analysis that the system $\bigl\{\epower{2\pi\ii \nu(t/(2a))}\,\big|\,\nu\in\Z\bigr\}$ is an orthonormal basis for the Hilbert space of square-integrable functions on $[-a, a]$ equipped with the standard inner product $\inprod{f}{g} \Let \frac{1}{2a}\int_{-a}^a f(t)g(t){\mrm dt}$. We let
			\begin{align*}
				\hat\HH & \Let \linspan\Biggl\{\biggl[\frac{\sin(\pi\nu t_1/a)}{\sqrt n}, \ldots, \frac{\sin(\pi\nu t_n/2)}{\sqrt n}\biggr]\transp\,\Bigg|\,t_i\in[-a,a], i=1, \ldots, n, \nu=1, \ldots, M\Biggr\}.
			\end{align*}
			Let $\ee^\nu(t_1, \ldots, t_n) \Let \sqrt{\frac{2}{n}}\bigl[\sin(\pi\nu t_1/a), \ldots, \sin(\pi\nu t_n/a)\bigr]\transp$, $t_i\in[-a, a]$. It is clear that the $\R^n$-valued functions $\bigl\{\ee^\nu,\;\nu=1, \ldots, M\bigr\}$ form an orthonormal set. Indeed,
			\begin{align*}
				\inprod{\ee_{\nu_1}}{\ee_{\nu_2}}_{\hat\HH} & = \sum_{i=1}^n\inprod{\ee_{\nu_1, i}}{\ee_{\nu_2, i}} = \frac{2}{n}\sum_{i=1}^n \frac{1}{2a}\int_{-a}^a\sin(\pi\nu_1 t_i/a)\sin(\pi\nu_2 t_i/a)\mrm dt_i\\
				& = \frac{2}{n}\sum_{i=1}^n\frac{1}{4}\int_{-1}^1 \bigl(\cos((\nu_1 - \nu_2)\pi s_i) - \cos((\nu_1 + \nu_2)\pi s_i)\bigr)\mrm ds_i\\
				& = \begin{cases}
					\frac{1}{2n}\sum_{i=1}^n 2 = 1 & \text{if }\nu_1 = \nu_2,\\
					0 & \text{otherwise}.
				\end{cases}
			\end{align*}
			We define $u_t \Let \eta_t + \Theta_t \ee(  w) = \eta_t + \sum_{j=0}^{t-1}\theta_{t, j} \ee(w_j) = \eta_t + \sum_{j=0}^{t-1}\sum_{\nu=1}^M\theta_{t, j}^\nu\ee^\nu(w_j)$ for appropriate matrices $\theta_{t, j}^\nu$. Now finding policies of the form~\eqref{e:hspolicy} that minimize the objective function in~\eqref{eq:problem1} becomes straightforward in the setting of Corollary~\ref{c:hspolicy}. The matrices $\Sigma_{\ee}$ and $\Sigma_{\ee}'$ in Corollary~\ref{c:hspolicy} are now easy to derive from Euler's identity $\epower{\ii\theta} = \cos\theta + \ii\sin \theta$, and the fact that the characteristic function of a uniform random variable $\zeta$ supported on $[-a, a]$ is given by $\EE\bigl[\epower{2\pi\ii v \zeta}\bigr] = \frac{1}{2a}\int_{-a}^a \epower{2\pi\ii v t}\,\mrm dt = \sinc(2\pi va)$ for some $v\in\R$, where the function $\sinc$ is defined as $\sinc(\xi) \Let \sin(\xi)/\xi$ if $\xi\neq 0$ and $1$ otherwise.

An alternative representation of the various matrices may be obtained by looking at each component of the policy elements separately. In this approach we define 
			\begin{align*}
				 {\hat\ee}(w_{t, i}) & \Let \begin{bmatrix} \ee_0(w_{t, i})& \ee_{1}(w_{t, i}) &  \cdots & \ee_{M}(w_{t, i}) \end{bmatrix}\transp,\\
				 {\hat\ee}(w_t) & \Let \begin{bmatrix} {\hat\ee}(w_{t, 1})\transp & \cdots & {\hat\ee}(w_{t, n})\transp\end{bmatrix}\transp, \quad  {\hat\ee}(  w) \Let \begin{bmatrix} {\hat\ee}(w_0)\transp & \cdots &{\hat\ee}(w_{N-2})\transp\end{bmatrix}\transp.
			\end{align*}
			In the above notation $\eta_{t, i} + \sum_{j=0}^{t-1}\theta_{j, i} {\hat\ee}(  w_{j, i})$ is of course the $i$-th entry of the input $u_t$ at time $t$, where $t=0, \ldots, N-1$ and $i=1, \ldots, n$.
\ExEnd
		\end{example}

	\subsection{Constraints on control energy}
		Some applications require constraints on the total control energy expended over a finite horizon. In the framework that we have established so far, such constraints are easy to incorporate. Indeed, if we require that $u\transp S u \le \beta^2$ for some preassigned $\beta > 0$ and positive definite matrix $S$, then in the setting of Theorem \ref{t:gen} this can be ensured by adjoining the condition $\norm{\eta}_{S} + \norm{\Theta}_{S}\norm{S}_\infty\mcal E \le \beta$ to the constraints, where $\norm{\eta}_{M} \Let \sqrt{\eta\transp M \eta}$ is the standard weighted $2$-norm for a positive definite matrix $M$.


\subsection*{Comparison with affine policies}

As pointed out earlier affine feedback policies from the noise have been previously treated in \cite{ref:loefberg03,ref:ben-tal04, ref:goulart06, GoulartKerrigan-08}, where the following feedback policy was considered:
		\begin{equation}
		\label{eq:bentalpolicy}
			u_t=\sum_{i=0}^{t-1} \theta_{t,i} w_i + \eta_t.
		\end{equation}
In the deterministic setting it was shown in~\cite{ref:goulart06} that there exists a one-to-one (nonlinear) mapping between control policies in the form~\eqref{eq:bentalpolicy} and the class of affine state feedback policies. That is, provided one is interested in affine state feedback policies, the parametrization~\eqref{eq:bentalpolicy} constitutes no loss of generality. In fact, we shall illustrate in the examples,   in the unconstrained inputs case, that the performance of this strategy with $\ee(w_i)$ in place of $w_i$ is almost as good as the standard LQG controller if not equally good.
However, in the constrained inputs case this choice is suboptimal in the class of measurable control policies, but it ensures tractability of a large class of optimal control problems.
		It can be seen that the solution to the optimization problem~\eqref{eq:problem} is tractable with this parametrization~\cite{ref:goulart06}. However, if the elements of the noise vector $ w$ are unbounded, the control input~\eqref{eq:bentalpolicy} does not have an upper bound. For the case of bounded inputs, the control policy~\eqref{eq:bentalpolicy} under unbounded noise will in general not satisfy the control constraint sets \eqref{eqn:Utotal}. This unboundedness is a potential problem in practical applications, and has been usually circumvented by assuming that the noise input lies within a compact set~\cite{BertsimasBrown-07, ref:goulart06} and designing a worst-case min-max type controller under this assumption.

		It is important to point out that our result in Section \ref{sec:bcbn} differs from that in \cite{ref:goulart06} in two aspects. First, we are solving the problem on finite-dimensional Hilbert spaces with general basis functions as opposed to a finite collection of affine functions in \cite{ref:goulart06}. Second, the feasibility of our problem is maintained for any bound on the elements of $\W$, as our constraint in \eqref{e:probhsrel} could still produce a feedback gain matrix $\Theta$ that has norm substantially different that $0$, whereas if there are elements in $\W$ with large enough norm and we take the control input to be $  u = \eta+ \Theta  w$, the constraints produce always a solution $\Theta$ with norm very close to $0$, hence practically only the open-loop term remains in the case of \cite{ref:goulart06}.


\section{Stability Analysis} \label{sec:stability}

	The main result in Theorem \ref{t:gen} asserts that the finite horizon optimization problem \eqref{eq:problem1} is convex and tractable using the policy \eqref{e:genpolicy}. To apply this result in a receding horizon fashion, it is imperative to further study some qualitative stability properties of the proposed policy. Under this policy the closed-loop system is not necessarily Markovian, and as such, standard Foster-Lyapunov methods cannot be directly applied. In what follows, we treat the stability problem for $p=\infty$ and $U_{\max}\Let U_{\max}^{(\infty)}$. However, this is without any loss of generality, for the same results hold (with minor modifications in the proofs) for $p = 1, 2$ as well. We impose the following assumption:
	\begin{assumption}
	\label{a:Schur}
		The matrix $A$ in~\eqref{eq:system} is Schur stable, i.e., the absolute value of the eigenvalues of $A$ are all strictly less than $1$.\AssmpEnd
	\end{assumption}
	At a first glance this assumption on $A$ might seem restrictive. Indeed, in the deterministic setting we know \cite{ref:YangSontagSussmann-97} that for discrete-time controlled systems it is possible to achieve global asymptotic stability with bounded control inputs if and only if the pair $(A, B)$ is stabilizable with arbitrary controls, and the spectral radius of $A$ is at most $1$. However, the problem of ensuring bounded variance of linear stochastic systems with bounded control inputs is to our knowledge still largely open; see, however, the recent manuscript \cite{ref:ramponinstable} for partial results as well as in \cite{ref:batina02, ref:Stoorvogel}.

	\subsection{Mean-square boundedness}
	We shall show that the variance of the state is uniformly bounded under receding horizon application of the strategy \eqref{e:genpolicy}, for any control horizon $N_c\le N$.
	The receding horizon implementation is iterative in nature: the optimization problem is solved every $kN_c$ steps, where $k\in\Nz$. The resulting optimal control policy (applied over a horizon $N_c$) is given by
	\[
		\pi_{kN_c:(k+1)N_c - 1}^*(x_{kN_c})\Let\smat{\pi^*_{kN_c}(x_{kN_c})\\\pi^*_{kN_c + 1}(x_{kN_c})\\\vdots\\\pi^*_{(k+1)N_c - 1}(x_{kN_c})} = \smat{\eta_0^*(x_{kN_c})\\\eta_1^*(x_{kN_c}) + \Theta_1^*(x_{kN_c})\ee(w)\\\vdots\\\eta_{N_c - 1}^*(x_{kN_c}) + \Theta^*_{N_c-1}(x_{kN_c})\ee(w)}
	\]
	where the control gains depend explicitly on the initial condition $x_{kN_c}$. For $\ell = 1,\cdots, N_c$, the resulting closed-loop system over horizon $N_c$ is given by:
\begin{equation}\label{eq:RHcl}
    x_{kN_c+\ell}=A^\ell x_{kN_c}+ B_\ell \pi_{kN_c:kN_c + \ell - 1}^*(x_{kN_c}) + D_\ell \tilde w_{kN_c:kN_c + \ell - 1},\qquad k\in\Nz,
\end{equation}
where $B_\ell \Let \smat{{\bar A}^{\ell-1}\bar B & \cdots & \bar A\bar B &\bar B}$, $D_\ell \Let \smat{{\bar A}^{\ell-1} & \cdots & \bar A & \mathbf I_{n\times n}}$, and $\tilde w_{kN_c:kN_c + \ell -1} \Let \smat{w_{kN_c}\transp & \cdots & w_{kN_c + \ell - 1}\transp}\transp$.

	Suppose that the above $N_c$-horizon optimal policy is computed as in Corollary \ref{c:sat}. We define the receding horizon policy corresponding to the consecutive concatenation of this $N_c$-horizon optimal policy as
	\begin{equation}
	\label{e:rhpolicy}
		\pi^* \Let \bigl(\pi_{0:N_c - 1}^*(x_{0}), \; \pi_{N_c:2N_c - 1}^*(x_{N_c}), \; \pi_{2N_c:3N_c - 1}^*(x_{2N_c}),  \cdots \bigr).
	\end{equation}

	\begin{proposition}\label{prop:RHmain}
		Consider the system~\eqref{eq:system}, and suppose that Assumptions~\ref{a:sys} and~\ref{a:Schur} hold. For $p = \infty$ and any control horizon $1\le N_c\le N$ the receding horizon control policy $\pi^*$ renders the closed loop system~\eqref{eq:RHcl} mean-square bounded, i.e., $\sup_{t\in\Nz} \EE_{\xz}\bigl[\norm{x_t}^2\bigr] < \infty$ for every initial condition $\xz\in\R^n$.
	\end{proposition}

	The proof of this Proposition is postponed to \secref{s:app:proof} in the Appendix.

\subsection{Input-to-state stability} \label{sec:iss}
	Input-to-state stability (\iss{}) is an interesting and important qualitative property of input-output behavior of dynamical systems. In the deterministic discrete-time setting \cite{JiangWang-01}, \iss{} generalizes the well-known bounded-input bounded-output (BIBO) property of linear systems~\cite[p.~490]{ref:antsaklisLS} to the setting of nonlinear systems. \iss{} provides a description of the behavior of a system subjected to bounded inputs, and as such it may be viewed as an $\mc L_\infty$ to $\mc L_\infty$ gain of a given nonlinear system. In this section we are interested in a useful stochastic variant of input-to-state stability; see e.g.,~\cite{ref:borkarUnifStab, ref:tsiniasRobustStochISS} for other possible definitions and ideas (primarily in continuous-time).


	\begin{defn}
		The system~\eqref{eq:system} is \emph{input-to-state stable in $\mathcal L_1$} if there exist functions $\beta\in\ClassKL$ and $\alpha, \gamma_1, \gamma_2\in\ClassKinfty$ such that for every initial condition $\xz\in\R^n$ we have
		\begin{equation}
		\label{eq:issm}
			\EE_{\xz}\bigl[\alpha(\norm{x_t})\bigr] \le \beta(\norm{\xz}, t) + \gamma_1\Bigl(\sup_{s\in\Nz}\norm{u_s}_\infty\Bigr) + \gamma_2\bigl(\norm{\Sigma}'\bigr)\qquad\fa t\in\Nz,
		\end{equation}
		where $\norm{\cdot}'$ is an appropriate matrix norm.\DefEnd
	\end{defn}

	One evident difference of \iss{} in $\mcal L_1$ with the deterministic definition of \iss{} is the presence of the function $\alpha$ \textsl{inside} the expectation in~\eqref{eq:issm}. It turns out that often it is more natural to arrive at an estimate of $\EE_{\xz}[\alpha(\norm{x_t})]$ for some $\alpha\in\ClassKinfty$ than an estimate of $\EE_{\xz}[\norm{x_t}]$. Moreover, in case $\alpha$ is convex, Jensen's inequality~\cite[p.~348]{ref:dudley} implies that such an estimate is stronger than an estimate of $\EE_{\xz}[\norm{x_t}]$.

	The property expressed by~\eqref{eq:issm} is one possible \iss{}-type property for stochastic systems. One can come up with alternative stochastic analogs of the \iss{} property, such as the following: $\fa\eps\in\;]0, 1[$ $\exists\,\beta\in\ClassKL$ and $\exists\,\gamma_1, \gamma_2\in\ClassKinfty$ such that $\PP\bigl(\norm{x_t} \le \beta(\norm{\xz}, t) + \gamma(\sup_{s\in\Nz}\norm{u_s}) + \gamma_2(\norm{\Sigma}')\;\fa t\in\Nz\bigr) \ge 1-\eps$. Intuitively this means that for $1-\eps$ proportion of the sample paths the deterministic \iss{} property holds uniformly. However, in an additive i.i.d unbounded noise setting as in~\eqref{eq:system}, this property fails to hold because almost surely the states undergo excursions outside any bounded set infinitely often; in this case the weaker version: $\fa\eps\in\;]0, 1[$ $\exists\,\beta\in\ClassKL$ and $\exists\,\gamma_1, \gamma_2\in\ClassKinfty$ such that $\PP\bigl(\norm{x_t} \le \beta(\norm{\xz}, t) + \gamma(\sup_{s\in\Nz}\norm{u_s}) + \gamma_2(\norm{\Sigma}')\bigr) \ge 1-\eps\;\fa t\in\Nz$ is comparatively better suited. We shall however stick with the \iss{} in $\mcal L_1$ property in this article.

	The following Proposition can be established with the aid of Proposition~\ref{prop:RHmain} for $p = \infty$; the proofs for $p = 1$ and $2$ are also similar in spirit.

	\begin{proposition}
		Consider the system~\eqref{eq:system}, and suppose that Assumptions~\ref{a:sys} and \ref{a:Schur} hold. Then the closed-loop system~\eqref{eq:RHcl} is \iss{} in $\mathcal L_1$ under the policy $\pi^*$ in \eqref{e:rhpolicy} for any $1 \le N_c \le N$.
	\end{proposition}

\section{Numerical Examples}\label{sec:examples}
	In this section we present several numerical examples to illustrate the theoretical results in the preceding sections. We start in Example \ref{ex:unconstrained} by comparing the performance of our policy \eqref{e:Thetadef} to that of the standard finite horizon LQG controller whenever the control inputs set $\bar \U\equiv\R^m$, i.e., there are no bounds on the norm of the inputs. Then we compare the performance of our policy \eqref{e:Thetadef} against a saturated LQG controller in Example \ref{ex:satLQGvsRHC}. Finally, in Example \ref{ex:constrained} we illustrate the effectiveness of our policy \eqref{e:Thetadef} compared to the certainty-equivalent receding horizon control.

	\begin{example}[Unconstrained Inputs]
	\label{ex:unconstrained}
		A natural question that may arise whenever the control inputs in our setup are not constrained, i.e., $\bar \U\equiv \R^m$, is the following: How does the policy \eqref{e:Thetadef} compare to the globally optimal controller, which in this case is the standard finite-horizon LQG controller? One would expect our policy to perform worse on the average since we restrict to a class of feedback policies that may not contain the globally optimal controller.

		We compared our policy against that of the LQG problem in simulation for two controllable $3$-dimensional single-input linear systems. In each case we solved an unconstrained finite-horizon LQ optimal control problem corresponding to state and control weights $Q_t = 3\,\mathbf I_{3\times 3}$ and $R_t = 1$ for every $t$.  We selected an optimization horizon $N=50$, and simulated the system responses starting from $10^3$ different initial conditions $x_0$ selected at random uniformly from the cube $[-100, 100]^3$, and noise sequences $w_t$ corresponding to i.i.d Gaussian noise of mean $0$ and (randomly chosen) variance
		\[
			\Sigma_w =
			\begin{bmatrix}
				2.830399255 & 5.491512606 & 3.612257417\\
				5.491512606 & 11.554870229 & 6.896706270\\
				3.612257417 & 6.896706270 & 4.625993264
			\end{bmatrix}.
		\]
		We selected the nonlinear bounded term $\ee( w)$ in our policy $u = \eta + \Theta \ee(w)$  to be a vector of scalar sigmoidal functions $\varphi(\xi) \Let 0.2\xi/\sqrt{1+0.04\xi^2}$ applied to each coordinate of the vector $w$. The covariance matrices $\Sigma_{\ee}$ and $\Sigma_{\ee'}$ that are required to solve the optimization problem \eqref{e:Thetadef} were computed empirically via classical Monte Carlo methods \cite[Section 3.2]{RobertCasella-04} using  $10^6$ i.i.d samples.

		The first system is described by:
		\begin{equation}
		\label{e:ExStablesys}
			x_{k+1} =
			\begin{bmatrix}
				0 & 1 & 0\\
				0 & 0 & 1\\
				0.4 & 0.5 & -0.25
			\end{bmatrix}x_k +
			\begin{bmatrix}
				0\\
				0\\
				1
			\end{bmatrix}u_k +w_k.
		\end{equation}
		The system pair $(A, B)$ is in Brunovsky canonical form, and $A$ has eigenvalues at $0.8642$, and $-0.5571 \pm\ii 0.3905$. The test results showed that the mean of the ratio of the cost corresponding to LQG to the cost corresponding to our policy is $0.99916$, and the standard deviation of this ratio is $0.003619$.


		The second system is described by:
		\begin{equation}
		\label{e:ExMarginalsys}
			x_{k+1}=
			\begin{bmatrix}
				1 & 1 & 0\\
				0 & 1 & 1\\
				0 & 0 & 1
			\end{bmatrix}x_k+
			\begin{bmatrix}
				0\\
				0\\
				1
			\end{bmatrix}u_k+w_k.
		\end{equation}
		This particular system matrix $A$ is in Jordan canonical form and has three eigenvalues at $1$. The test results showed that the mean of the ratio of the cost of LQG against the cost of our policy is $0.99673$ and the corresponding standard deviation is $0.008045$.

		Computations for determining our policy in the above two cases were carried out in the MATLAB-based software package \texttt{cvx}. In the case of the system \eqref{e:ExMarginalsys} the solver utilized by \texttt{cvx} reported numerical problems in five different runs, for which it gave values of the aforementioned ratio below $0.96$. Note that we have not discarded these five cases from the mean and variance figures reported above.

		The close-to-optimal performance of our policy is surprising in view of the fact that the vector-space $\HH$ is the linear span of one bounded function, and does not contain the theoretically optimal linear (in the current state) controller. We conjecture that this is due to injectivity of the mapping $\ee$, due to which $\ee(w_t)$ retains all information generated by $w_t$. Of course, in the absence of control constraints our solution is much more computationally demanding than the LQG controller, and would not be used in practice in this case.\ExEnd
	\end{example}

	\begin{example}[Saturated LQG and Receding Horizon]
	\label{ex:satLQGvsRHC}
		We compare the performance of saturated LQG against our policy \eqref{e:Thetadef} for the system \eqref{e:ExMarginalsys} in this example. We fixed the optimization horizon $N=2$, the control horizon $N_c=1$, and the weight matrices for the states and the control to be $Q_t=\mathbf I_{3\times 3}$ and $R_t=0.01$ for all $t$, respectively. The control bounds in both cases was $[-2, 2]$, the nonlinear bounded term $\ee(w_t)$ in our policy $u = \eta + \Theta\ee(w)$ was a vector of scalar standard saturation functions applied to each coordinate of the vector $w_t$, and the LQG control input was saturated at $\pm 2$. The covariance matrices $\Sigma_{\ee}$ and $\Sigma_{\ee'}$ required to solve the optimization problem \eqref{e:probspan} were computed empirically via classical Monte Carlo integration methods \cite[Section 3.2]{RobertCasella-04} using $10^6$ i.i.d samples.

		We simulated the system \eqref{e:ExMarginalsys} starting from the same initial condition $x_0=\smat{0& 0 & 0}\transp$ for $100$ different independent realizations of the noise sequence $w_t$ over a horizon of $200$. The behavior of the average (over the $100$ realizations) cost corresponding to the two scenarios is shown in Figure \ref{fig:LQGAffine}. The simulations were coded in MATLAB and the optimization programs were coded in the software package \texttt{cvx}. The average total cost incurred at the end of the simulation horizon when using the saturated LQG scheme above was $1.790\times 10^{12}$ units, whereas the average total cost incurred at the end of the simulation horizon ($t=200$) using our policy \eqref{e:Thetadef} in a receding horizon fashion was $4.486\times 10^8$ units.\ExEnd
		\begin{figure}[h]
			\centering
\includegraphics[width=0.75\columnwidth]{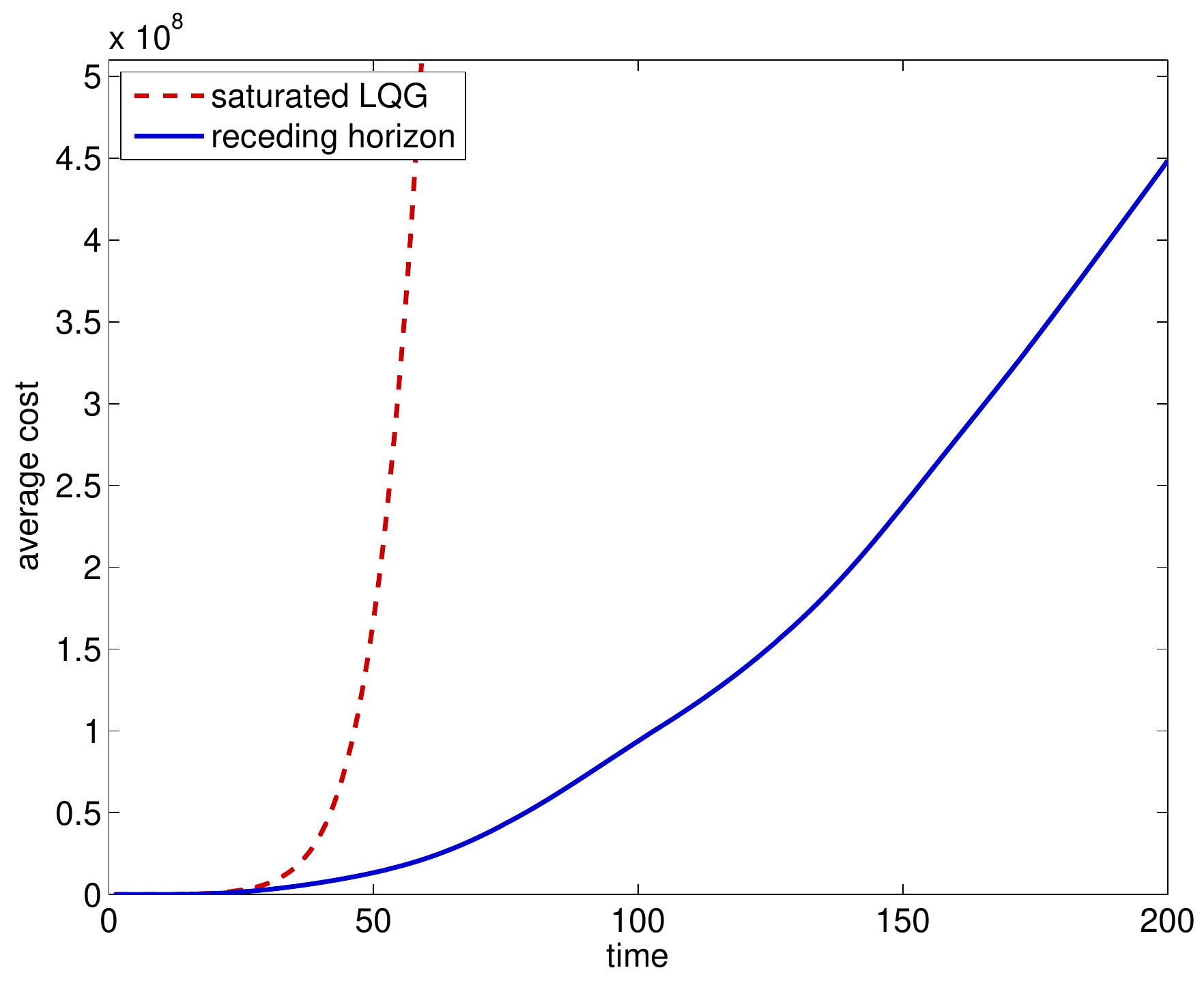}
			\caption{Plots of average costs corresponding to saturated LQG and our receding horizon scheme for $N_c = 1$ in Example \ref{ex:satLQGvsRHC}.}
			\label{fig:LQGAffine}
		\end{figure}
	\end{example}

	\begin{example}[Constrained Inputs]
	\label{ex:constrained}
		Consider the 2-dimensional linear stochastic system:
		\begin{equation}
		\label{e:crazysys}
			x_{t+1} =
			\begin{bmatrix}
				1.23	& -0.15\\
				0.25	& 1
			\end{bmatrix}
			x_t +
			\begin{bmatrix}
				0.14\\
				0.12
			\end{bmatrix}
			u_t + w_t,
		\end{equation}
		where $(w_t)_{t\in\Nz}$ is a sequence of i.i.d Gaussian noise with zero mean and (randomly generated) variance $\begin{bmatrix}2.722030613  &  4.975999693\\ 4.975999693  &  9.102559685\end{bmatrix}$. Let the weight matrices corresponding to the states and control be $Q_t = \mathbf I_{2\times 2}$ and $R_t = 0.8$ for each $t$. The covariance matrices $\Sigma_{\ee}$ and $\Sigma_{\ee'}$ that are required to solve the optimization problem \eqref{e:Thetadef} were computed empirically via classical Monte Carlo integration methods \cite[Section 3.2]{RobertCasella-04} using  $10^6$ samples.

		We fixed the optimization horizon $N = 7$, the nonlinear saturation $\ee(w_t)$ to be a vector of scalar sigmoidal functions $\varphi(\xi) \Let 0.2\xi/\sqrt{1+0.04\xi^2}$ applied to each coordinate of the vector $w_t$, and compared the certainty-equivalent MPC strategy ($N_c=1$, $\Theta\equiv 0$, $w_t\equiv 0$) against our receding horizon strategy \eqref{e:Thetadef} with control horizon $N_c=4$. The control constraints in both cases were $u_t \in [-200, 200]$. We simulated the system in both cases starting  from the same initial condition $x_0=\smat{0& 0}\transp$, for $60$ different realizations of the noise sequence $w_t$; plots of states, average cost, and standard deviation are shown in Figures \ref{fig:plotsmpc1rh4} and \ref{fig:plotsmpc1rh4costs}. The average cost incurred when using the certainty-equivalent MPC scheme was $7.893\times 10^5$ units, whereas the average cost incurred when using our policy \eqref{e:Thetadef} in a receding horizon fashion was $3.141\times 10^5$ units. Therefore, applying our policy in a receding horizon fashion one saves $60.2\%$ of the cost corresponding to the certainty-equivalent MPC controller on the average. This example illustrates that there may be cases where open-loop certainty-equivalent MPC, in the absence of state-constraints, is outperformed by a large margin by a judiciously selected receding-horizon strategy. The simulations were coded in \texttt{YALMIP} and were solved using \texttt{SDPT-3}; the solver-time statistics (in sec.) for the certainty-equivalent MPC and receding horizon schemes were as follows:\\
		\begin{center}
		\begin{tabular}{r|c|c}
			& certainty-equivalent MPC & receding horizon\\\hline
			Mean & $32.127$ & $59.615$\\
			Standard deviation & $4.610$ & $21.675$\\
			Maximum & $50.590$ & $90.036$\\
			Minimum & $20.240$ & $20.466$\\
		\end{tabular}
		\end{center}

\noindent $\phantom{s}$\\ These statistics correspond to the above simulations carried out on an $\text{x}86\_64$ octa-core machine with 24GB RAM, each processor of which was an Intel$^\text{\textregistered}$ Xeon$^\text{\textregistered}$ CPU E5540 \@ 2.53GHz with cache size 8192 KB, running GNU/Linux.
		\begin{figure}[h]
			\centering
			\includegraphics[width=0.8\columnwidth]{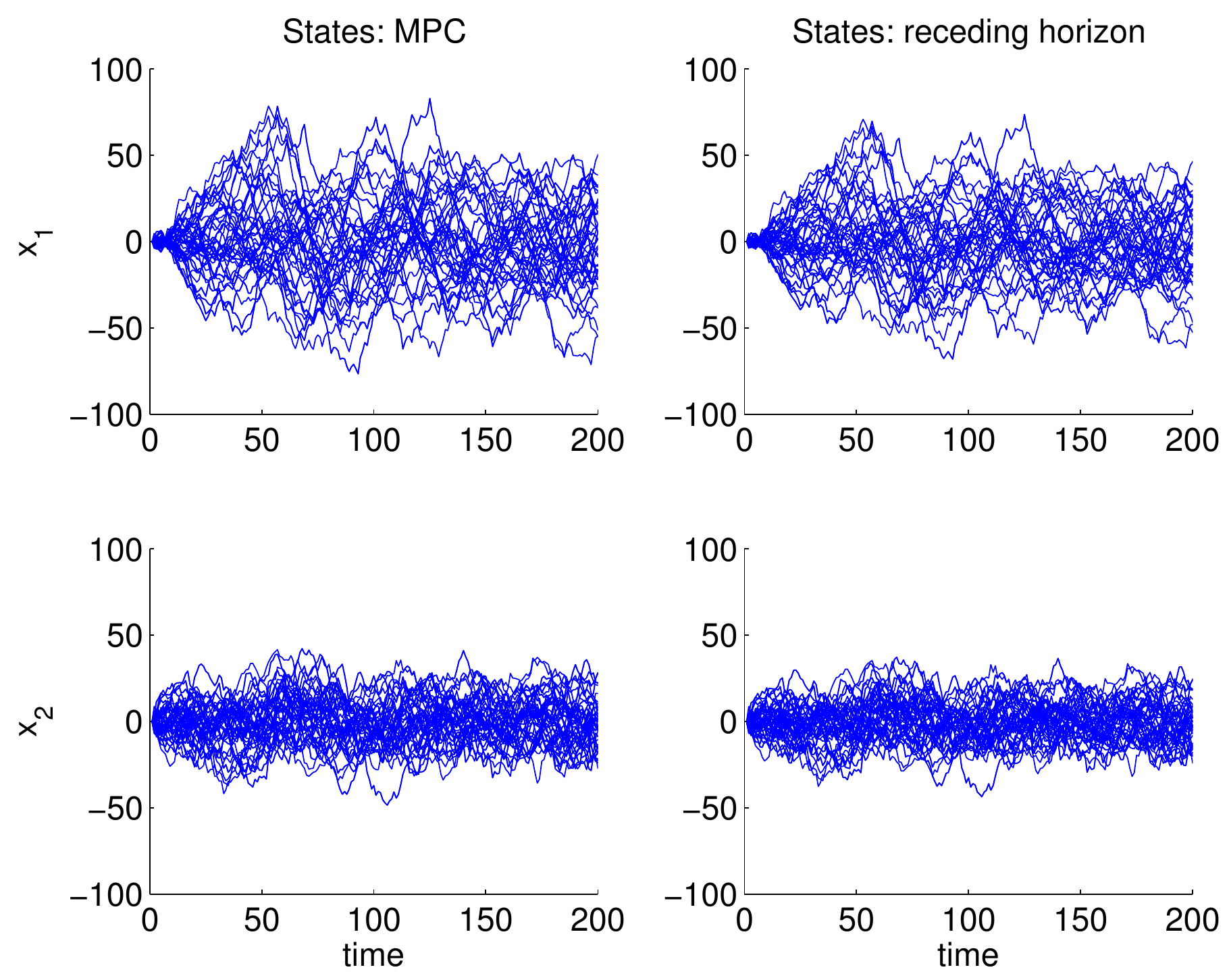}
			\caption{Plots of states corresponding to: certainty-equivalent MPC with $N_c=1$ (left) and our receding horizon control scheme with $N_c=4$ (right) in Example \ref{ex:constrained}.}
			\label{fig:plotsmpc1rh4}
			\label{plotsmpc1rh4}
		\end{figure}
		\begin{figure}[h]
			\centering
			\subfigure[Plot of average costs]{\includegraphics[width=0.48\columnwidth]{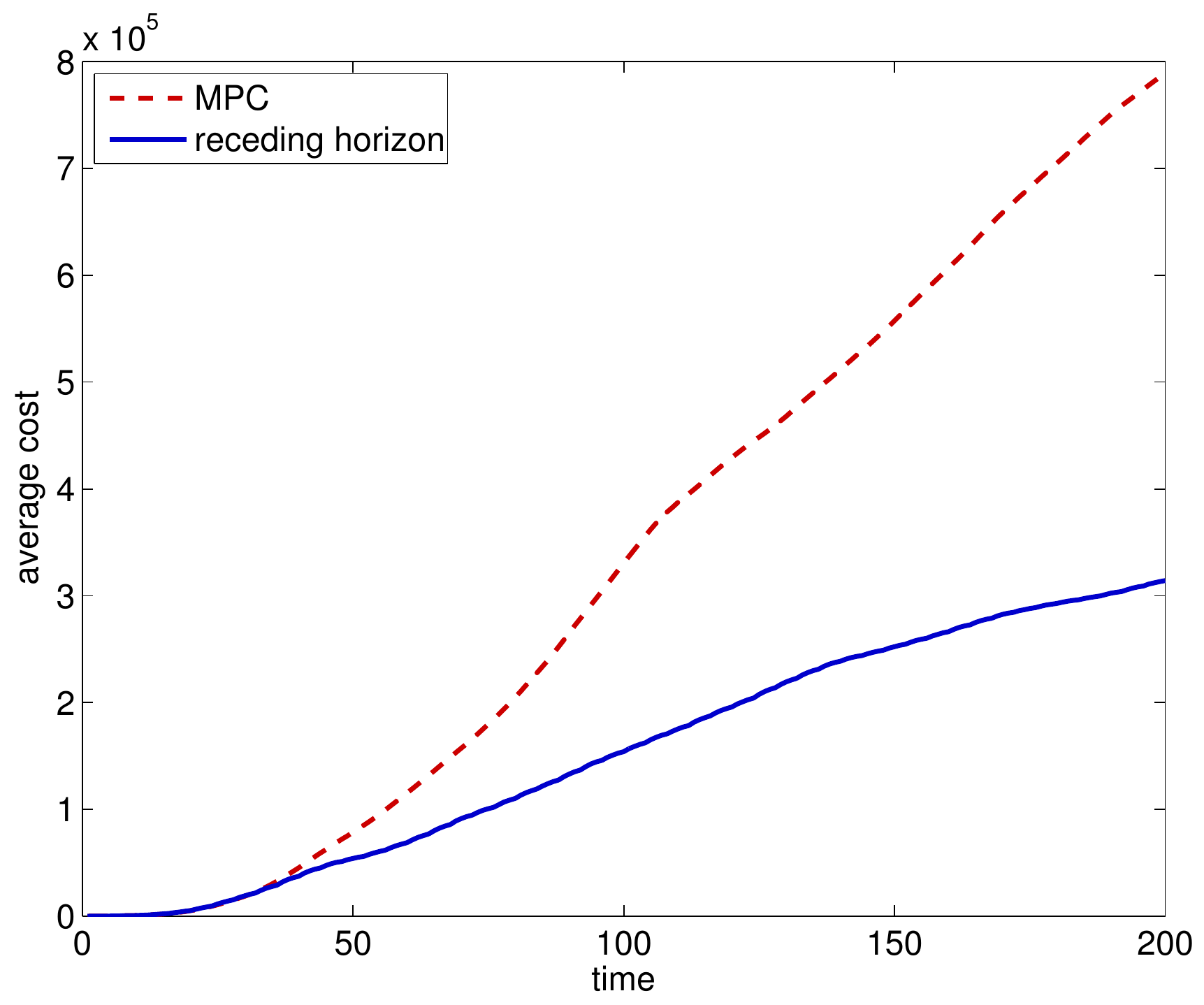}}
			\subfigure[Plot of standard deviations]{\includegraphics[width=0.48\columnwidth]{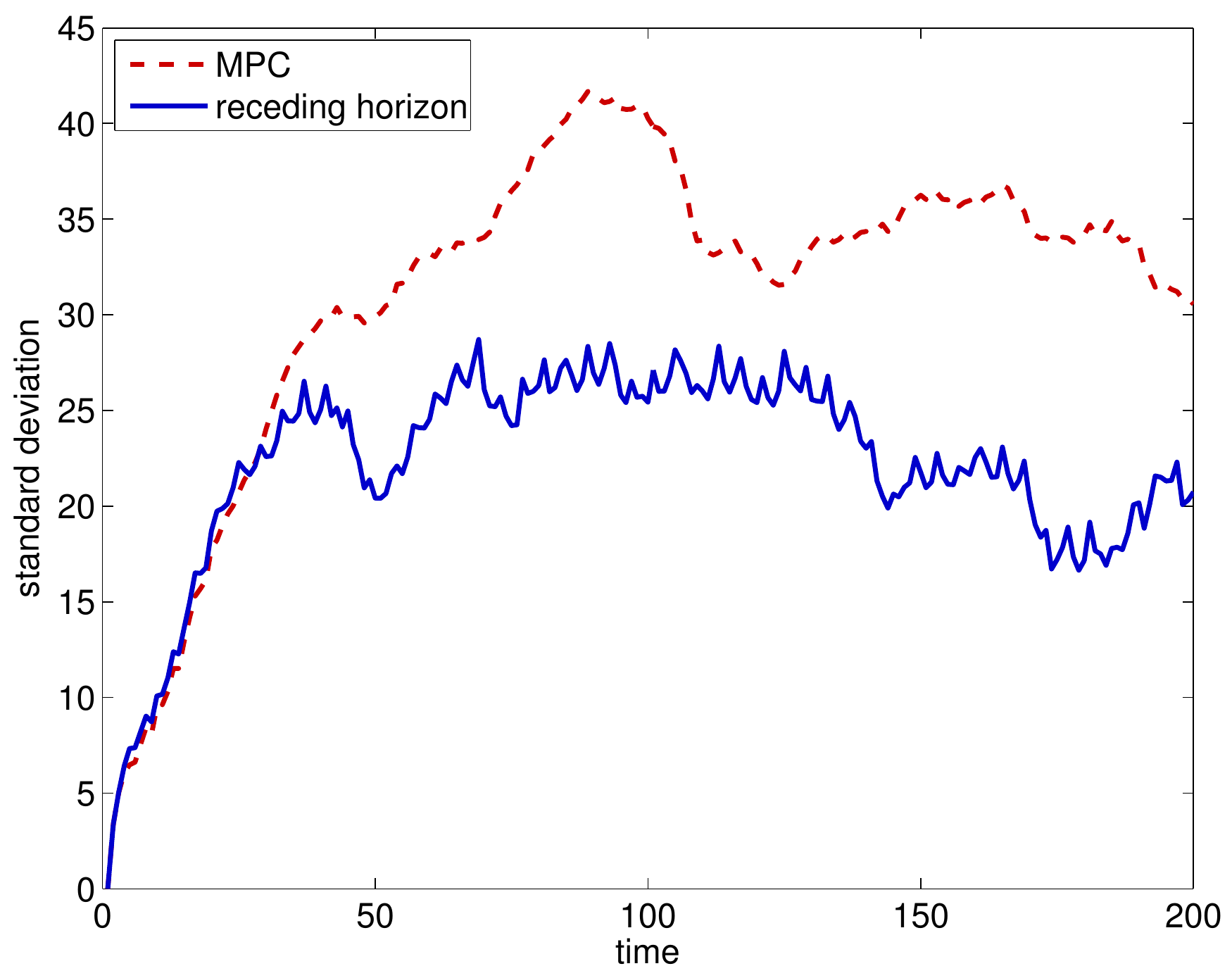}}
			\caption{Plots of average cost (left) and standard deviations (right) corresponding to: certainty-equivalent MPC with $N_c=1$ and our receding horizon control scheme with $N_c=4$ in Example \ref{ex:constrained}.}
			\label{fig:plotsmpc1rh4costs}
			\label{plotsmpc1rh4costs}
		\end{figure}

		We also applied the first four control values of the certainty-equivalent scheme and compared it against our receding horizon scheme using policy \eqref{e:Thetadef}, i.e., $N_c = 4$ for both controllers. We simulated the system in both cases starting  from the same initial condition $x_0=\smat{0 & 0}\transp$, for $60$ different realizations of the noise sequence $w_t$; plots of the states, average cost, and standard deviation are shown in Figures \ref{fig:plotsmpc4rh4} and \ref{fig:plotsmpc4rh4costs}. The average cost incurred when using the certainty-equivalent with control horizon $N_c=4$ was $4.211\times 10^5$ units, whereas the average cost incurred when using our policy \eqref{e:Thetadef} in a receding horizon fashion was $3.295\times 10^5$ units. We see that by applying our policy in a receding horizon fashion one saves $21.7\%$ of the cost corresponding to the certainty equivalence controller on the average. The simulations were coded in \texttt{YALMIP} and were solved using \texttt{SDPT-3}; the solver-time statistics (in sec.) for the certainty-equivalent and receding horizon schemes were as follows:\\
\begin{center}
		\begin{tabular}{r|c|c}
			& certainty-equivalent & receding horizon\\\hline
			Mean & $7.537$ & $67.494$\\
			Standard deviation & $0.812$ & $11.845$\\
			Maximum & $9.776$ & $85.232$\\
			Minimum & $6.101$ & $43.601$\\ 
		\end{tabular}
		\end{center}
		
\noindent $\phantom{s}$\\ These statistics correspond to the above simulations carried out on an $\text{x}86\_64$ octa-core machine with 24GB RAM, each processor of which was an Intel$^\text{\textregistered}$ Xeon$^\text{\textregistered}$ CPU E5540 \@ 2.53GHz with cache size 8192 KB, running GNU/Linux.\ExEnd

		\begin{figure}[h]
			\centering
			\includegraphics[width=0.8\columnwidth]{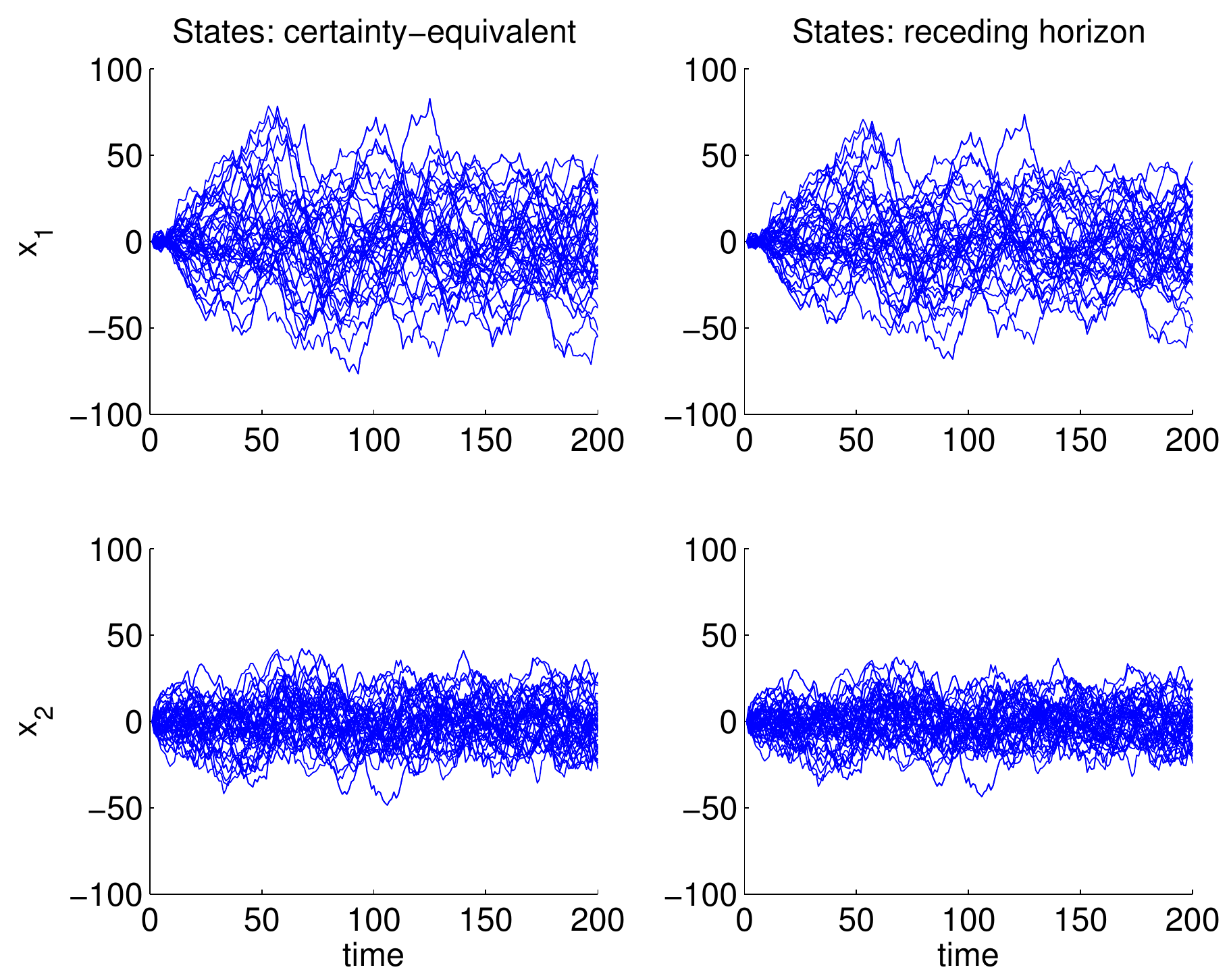}
			\caption{Plots of states corresponding to: certainty-equivalent with $N_c=4$ (left) and our receding horizon control scheme with $N_c=4$ (right) in Example \ref{ex:constrained}.}
			\label{fig:plotsmpc4rh4}
			\label{plotsmpc4rh4}
		\end{figure}
		\begin{figure}[h]
			\centering
			\subfigure[Plot of average costs]{\includegraphics[width=0.48\columnwidth]{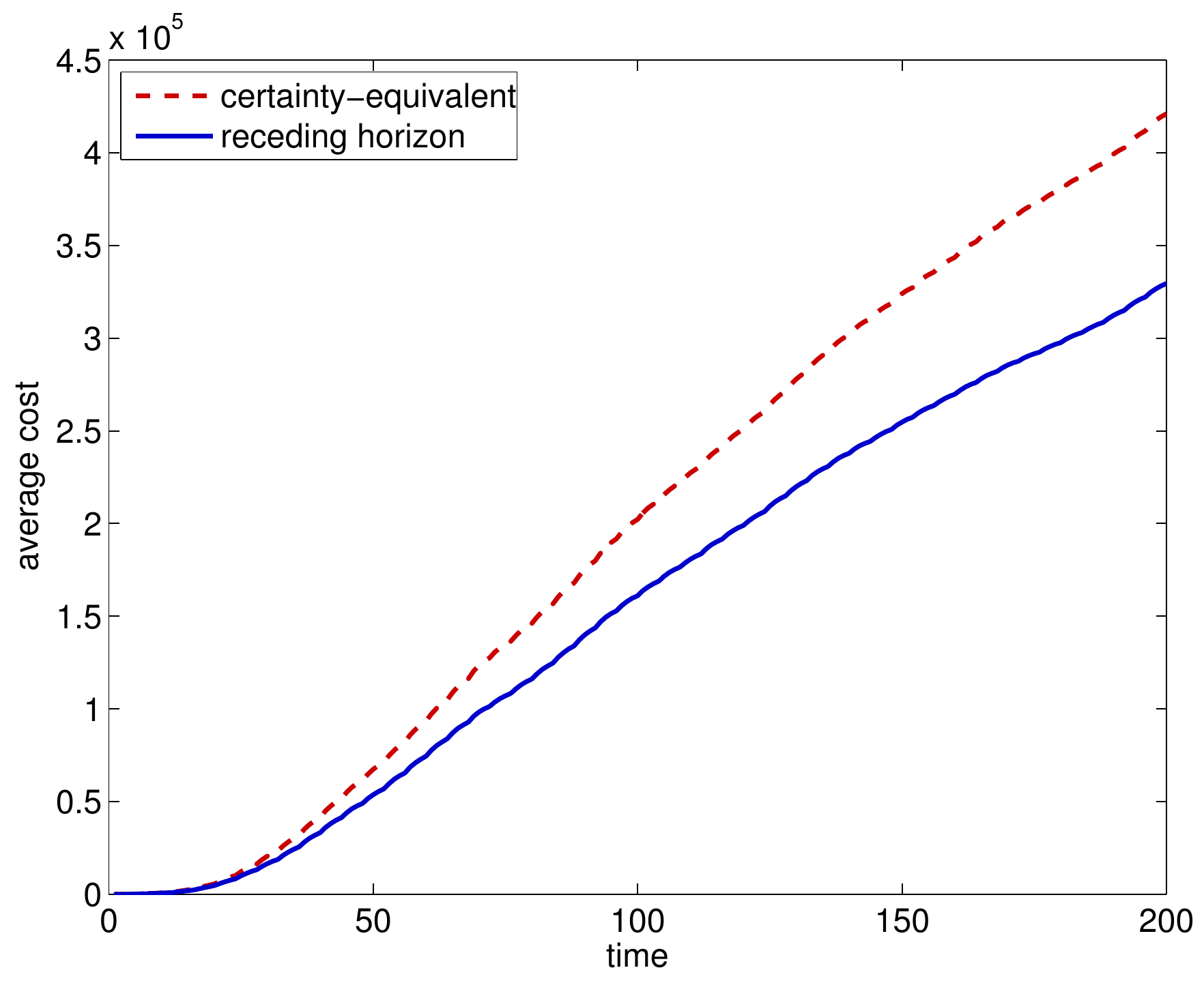}}
			\subfigure[Plot of standard deviations]{\includegraphics[width=0.48\columnwidth]{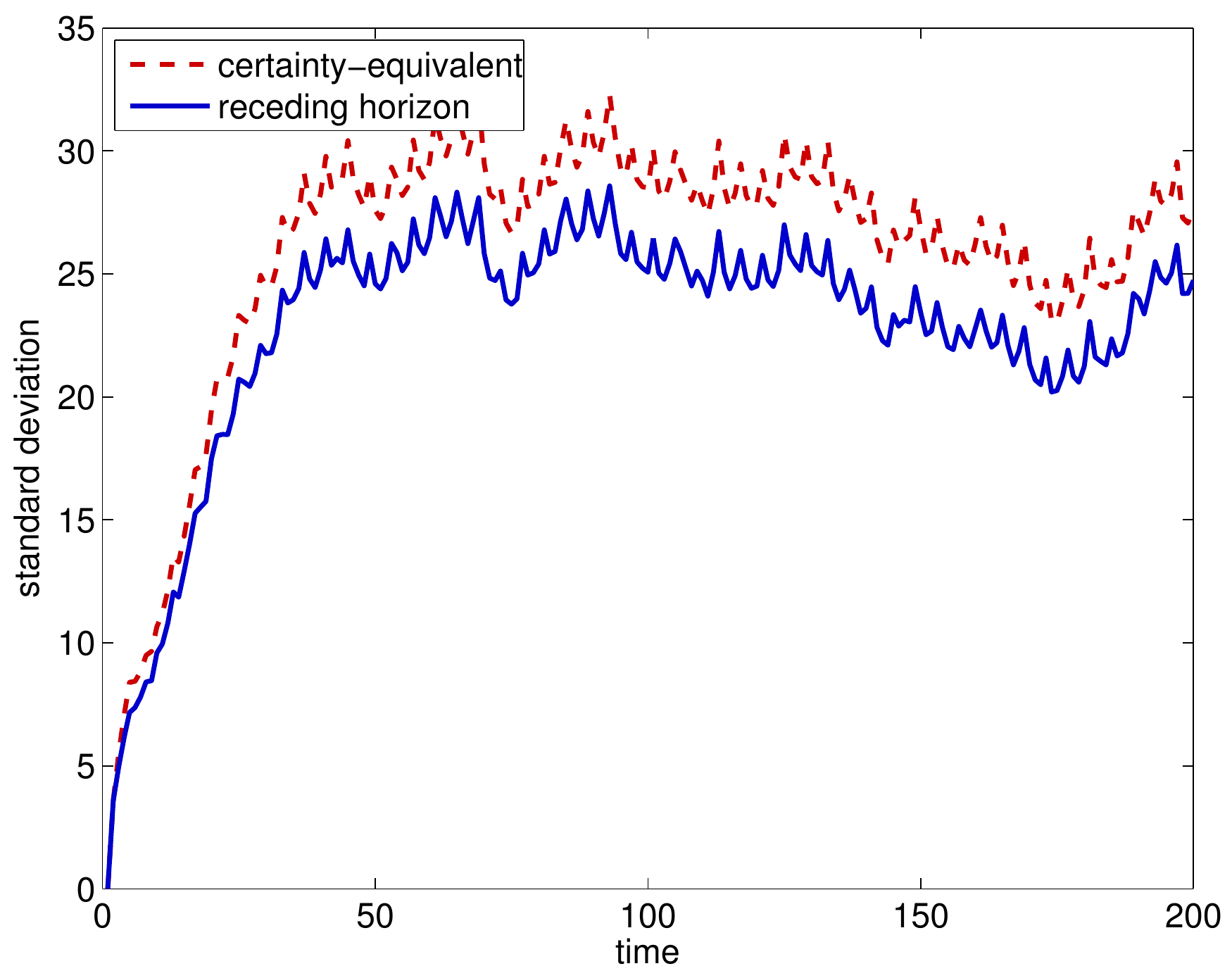}}
			\caption{Plots of average cost (left) and standard deviations (right) corresponding to: certainty-equivalent with $N_c=4$ and our receding horizon control scheme with $N_c=4$ in Example \ref{ex:constrained}.}
			\label{fig:plotsmpc4rh4costs}
			\label{plotsmpc4rh4costs}
		\end{figure}
\end{example}


\section{Conclusion and Future Directions}
\label{sec:conclusions}
We provided tractable solutions to a variety of finite-horizon stochastic optimal control problems with quadratic cost, hard control constraints, and unbounded additive noise. These problems arise as parts of solutions to the stochastic receding horizon problems~\eqref{eq:problem}. The control policy obtained as a result of the finite-horizon optimal control sub-problems may be nonlinear with respect to the previous states, and the policy elements are chosen from a vector space that is largely up to the designer. One of the key features of our approach is that the variance-like matrices employed in the finite-horizon optimal control sub-problems may be computed off-line, and we illustrated this feature with several examples. We demonstrated that applying our obtained policies in a receding horizon fashion results in bounded state variance. Finally, we provided several numerical examples that illustrate the effectiveness of our method with respect to the commonly used certainty-equivalent MPC controllers.

	The development in this article affords extensions in several directions. One is the incorporation of state constraints. As discussed in \secref{s:intro}, hard state constraints do not make sense in the stochastic with additive unbounded noise setting unless one is prepared to artificially relax them once infeasibility is encountered. Probabilistic constraints and integrated chance constraints~\cite{ref:haneveldICC} constitute popular alternative methods to impose constraints on the state that are more probabilistic in nature. It will be interesting to see how the approach introduced in this article reacts to state-constraints. A second direction is to consider specific kinds of nonlinear models, particularly those which involve multiplicative noise, in our framework, and a third is to consider different objective functions such as affine functions given by the $\ell_\infty$ and the $\ell_1$ norms.


\section*{Acknowledgments}
	We are indebted to Soumik Pal for pointing out the possibility of representing policies as elements of a vector space. We thank Colin Jones for some useful discussions on convexity of some of the optimization programs, and the three anonymous reviewers for their valuable suggestions that have led to substantial improvements of the original manuscript.


\begin{thebibliography}{RCMA{\etalchar{+}}09}

\bibitem[ACCL09]{AgarwalCinquemaniChatterjeeLygeros-09}
M.~Agarwal, E.~Cinquemani, D.~Chatterjee, and J.~Lygeros, \emph{On convexity of
  stochastic optimization problems with constraints}, European Control
  Conference, 2009, pp.~2827--2832.

\bibitem[AM06]{ref:antsaklisLS}
P.~J. Antsaklis and A.~N. Michel, \emph{Linear {S}ystems}, Birkh\"auser Boston
  Inc., Boston, MA, 2006.

\bibitem[AS64]{ref:AbramowitzStegun}
M.~Abramowitz and I.~A. Stegun, \emph{Handbook of {M}athematical {F}unctions
  with {F}ormulas, {G}raphs, and {M}athematical {T}ables}, National Bureau of
  Standards Applied Mathematics Series, vol.~55, Superintendent of Documents,
  U.S.\ Government Printing Office, Washington D.C., 1964.

\bibitem[Bat04]{batinaPhDthesis}
I.~Batina, \emph{Model predictive control for stochastic systems by randomized
  algorithms}, Ph.D. thesis, Technische Universiteit Eindhoven, 2004.

\bibitem[BB07]{BertsimasBrown-07}
D.~Bertsimas and D.~B. Brown, \emph{Constrained stochastic {LQC}: a tractable
  approach}, IEEE Transactions on Automatic Control \textbf{52} (2007), no.~10,
  1826--1841.

\bibitem[Ber05]{ref:bertsekas05survey}
D.~P. Bertsekas, \emph{Dynamic programming and suboptimal control: {a} survey
  from {ADP} to {MPC}}, European Journal of Control \textbf{11} (2005),
  no.~4-5, 310--334.

\bibitem[Ber09]{ref:bernsteinMatrixMath}
D.~S. Bernstein, \emph{Matrix {M}athematics}, 2 ed., Princeton University
  Press, 2009.

\bibitem[Bla99]{ref:blanchini1999sic}
F.~Blanchini, \emph{{Set invariance in control}}, Automatica \textbf{35}
  (1999), no.~11, 1747--1767.

\bibitem[BM99]{BemporadMorari-99}
A.~Bemporad and M.~Morari, \emph{Robust model predictive control: a survey},
  Robustness in Identification and Control \textbf{245} (1999), 207--226.

\bibitem[Bor00]{ref:borkarUnifStab}
V.~S. Borkar, \emph{Uniform stability of controlled {M}arkov processes}, System
  theory: modeling, analysis and control (Cambridge, MA, 1999), Kluwer
  International Series in Engineering Computer Science, vol. 518, Kluwer
  Academic Publishers, Boston, MA, 2000, pp.~107--120.

\bibitem[Bro97]{ref:brockettMinAttention}
R.~W. Brockett, \emph{Minimum attention control}, Proceedings of the 36th IEEE
  Conference on Decision and Control, vol.~3, 1997, pp.~2628--2632.

\bibitem[BSM03]{ref:bermanCPM}
A.~Berman and N.~Shaked-Monderer, \emph{Completely {P}ositive {M}atrices},
  World Scientific Publishing Co. Inc., River Edge, NJ, 2003.

\bibitem[BSW02]{ref:batina02}
I.~Batina, A.~A. Stoorvogel, and S.~Weiland, \emph{Optimal control of linear,
  stochastic systems with state and input constraints}, Proceedings of the 41st
  IEEE Conference on Decision and Control, vol.~2, 2002, pp.~1564--1569.

\bibitem[BT96]{ref:bertsekasNDP}
D.~Bertsekas and J.~Tsitsiklis, \emph{Neuro-{D}ynamic {P}rogramming}, Athena
  Scientific, 1996.

\bibitem[BTGGN04]{ref:ben-tal04}
A.~Ben-Tal, A.~Goryashko, E.~Guslitzer, and A.~Nemirovski, \emph{Adjustable
  robust solutions of uncertain linear programs}, Mathematical Programming
  \textbf{99} (2004), no.~2, 351--376.

\bibitem[BV04]{ref:boyd04}
S.~Boyd and L.~Vandenberghe, \emph{Convex {O}ptimization}, Cambridge University
  Press, Cambridge, 2004, Sixth printing with corrections, 2008.

\bibitem[BW07]{BlackmoreWilliams-07}
L.~Blackmore and B.~C. Williams, \emph{Optimal, robust predictive control of
  nonlinear systems under probabilistic uncertainty using particles},
  Proceedings of the American Control Conference, 2007, pp.~1759--1761.

\bibitem[CACL09]{ref:cinquemaniagarwal}
E.~Cinquemani, M.~Agarwal, D.~Chatterjee, and J.~Lygeros, \emph{On convex
  problems in chance-constrained stochastic model predictive control},
  \url{http://arxiv.org/abs/0905.3447}, 2009.

\bibitem[CCCL08]{ref:recstrat}
D.~Chatterjee, E.~Cinquemani, G.~Chaloulos, and J.~Lygeros, \emph{Stochastic
  optimal control up to a hitting time: optimality and rolling-horizon
  implementation}, \url{http://arxiv.org/abs/0806.3008}, 2008.

\bibitem[CKW08]{ref:CannonKouvaritakisWu-08}
M.~Cannon, B.~Kouvaritakis, and X.~Wu, \emph{Probabilistic constrained {MPC}
  for systems with multiplicative and additive stochastic uncertainty}, IFAC
  World Congress (Seoul, Korea), 2008.

\bibitem[CP09]{ref:excursion}
D.~Chatterjee and S.~Pal, \emph{An excursion-theoretic view of stability of
  stochastic hybrid systems}, \url{http://arxiv.org/abs/0901.2269}, 2009.

\bibitem[dFR03]{ref:vanRoyLPtoDP}
D.~P. de~Farias and B.~Van Roy, \emph{The linear programming approach to
  approximate dynamic programming}, Operations Research \textbf{51} (2003),
  no.~6, 850--865.

\bibitem[Dud02]{ref:dudley}
R.~M. Dudley, \emph{Real {A}nalysis and {P}robability}, Cambridge Studies in
  Advanced Mathematics, vol.~74, Cambridge University Press, Cambridge, 2002,
  Revised reprint of the 1989 original.

\bibitem[FB05]{FukushimaBitmead-05}
H.~Fukushima and R.~R. Bitmead, \emph{Robust constrained predictive control
  using comparison model}, Automatica \textbf{41} (2005), no.~1, 97--106.

\bibitem[GB00]{ref:boydCVX}
M.~Grant and S.~Boyd, \emph{{CVX}: {M}atlab software for disciplined convex
  programming (web page and software)}, \url{http://stanford.edu/~boyd/cvx},
  2000.

\bibitem[GK08]{GoulartKerrigan-08}
P.~J. Goulart and E.~C. Kerrigan, \emph{Input-to-state stability of robust
  receding horizon control with an expected value cost}, Automatica \textbf{44}
  (2008), no.~4, 1171--1174.

\bibitem[GKM06]{ref:goulart06}
P.~J. Goulart, E.~C. Kerrigan, and J.~M. Maciejowski, \emph{Optimization over
  state feedback policies for robust control with constraints}, Automatica
  \textbf{42} (2006), no.~4, 523--533.

\bibitem[Han83]{ref:haneveldICC}
W.~K.~Klein Haneveld, \emph{On integrated chance constraints}, Stochastic
  programming (Gargnano), Lecture Notes in Control and Inform. Sci., vol.~76,
  Springer, Berlin, 1983, pp.~194--209.

\bibitem[HCCL10]{ref:acc10}
P.~Hokayem, E.~Cinquemani, D.~Chatterjee, and J.~Lygeros, \emph{Stochastic
  {MPC} with output feedback and bounded control inputs}, 2010, Submitted to
  the American Control Conference.

\bibitem[HCL09]{ref:petercdc09}
P.~Hokayem, D.~Chatterjee, and J.~Lygeros, \emph{On stochastic model predictive
  control with bounded control inputs}, \url{http://arxiv.org/abs/0902.3944},
  2009.

\bibitem[JW01]{JiangWang-01}
Z.-P. Jiang and Y.~Wang, \emph{Input-to-state stability for discrete-time
  nonlinear systems}, Automatica \textbf{37} (2001), no.~6, 857--869.

\bibitem[LH07]{lavretsky2007stable}
E.~Lavretsky and N.~Hovakimyan, \emph{{Stable adaptation in the presence of
  actuator constraints with flight control applications}}, Journal of Guidance
  Control and Dynamics \textbf{30} (2007), no.~2, 337.

\bibitem[LHBW07]{ref:lazarbemporad07}
M.~Lazar, W.~P. M.~H. Heemels, A.~Bemporad, and S.~Weiland, \emph{Discrete-time
  non-smooth nonlinear {MPC}: stability and robustness}, Lecture Notes in
  Control and Information Sciences, vol. 358, Springer-Verlag, 2007,
  pp.~93--103.

\bibitem[LHC03]{LavretskyHovakimyanCalise-03}
E.~Lavretsky, N.~Hovakimyan, and A.~J. Calise, \emph{Upper bounds for
  approximation of continuous-time dynamics using delayed outputs and
  feedforward neural networks}, IEEE Transactions on Automatic Control
  \textbf{48} (2003), no.~9, 1606--1610.

\bibitem[L{\"o}f03]{ref:loefberg03}
J.~L{\"o}fberg, \emph{Minimax {A}pproaches to {R}obust {M}odel {P}redictive
  {C}ontrol}, Ph.D. thesis, Link{\"o}pings Universitet, 2003.

\bibitem[L{\"o}f04]{YALMIP}
\bysame, \emph{{YALMIP} : {A} {T}oolbox for {M}odeling and {O}ptimization in
  {MATLAB}}, Proceedings of the CACSD Conference (Taipei, Taiwan), 2004.

\bibitem[LR06]{ref:rantzerRelaxingDP}
B.~Lincoln and A.~Rantzer, \emph{Relaxing dynamic programming}, IEEE
  Transactions on Automatic Control \textbf{51} (2006), no.~8, 1249--1260.

\bibitem[Lue69]{Luenberger-69}
D.~G. Luenberger, \emph{Optimization by {V}ector {S}pace {M}ethods}, J. Wiley
  \& Sons, 1969.

\bibitem[Mac01]{ref:maciejowskibk}
J.~M. Maciejowski, \emph{Predictive {C}ontrol with {C}onstraints}, Prentice
  Hall, 2001.

\bibitem[MLL05]{MaciejowskiLecchiniLygeros-05}
J.~M. Maciejowski, A.~Lecchini, and J.~Lygeros, \emph{{NMPC} for complex
  stochastic systems using {Markov Chain Monte Carlo}}, International Workshop
  on Assessment and Future Directions of Nonlinear Model Predictive Control
  (Stuttgart, Germany), Lecture Notes in Control and Information Sciences, vol.
  358/2007, Springer, 2005, pp.~269--281.

\bibitem[MRRS00]{MayneRawlingsRaoScokaert-00}
D.~Q. Mayne, J.~B. Rawlings, C.~V. Rao, and P.~O.~M. Scokaert,
  \emph{Constrained model predictive control: stability and optimality},
  Automatica \textbf{36} (2000), no.~6, 789--814.

\bibitem[OJM08]{OldewurtelJonesMorari-08}
F.~Oldewurtel, C.N. Jones, and M.~Morari, \emph{A tractable approximation of
  chance constrained stochastic {MPC} based on affine disturbance feedback},
  Proceedings of the 47th IEEE Conference on Decision and Control, 2008,
  pp.~4731--4736.

\bibitem[Pow07]{ref:powellADP}
W.~B. Powell, \emph{Approximate {D}ynamic {P}rogramming}, Wiley Series in
  Probability and Statistics, Wiley-Interscience [John Wiley \& Sons], Hoboken,
  NJ, 2007.

\bibitem[PS09]{PrimbsSung-09}
J.~A. Primbs and C.~H. Sung, \emph{Stochastic receding horizon control of
  constrained linear systems with state and control multiplicative noise}, IEEE
  Transactions on Automatic Control \textbf{54} (2009), no.~2, 221--230.

\bibitem[RC04]{RobertCasella-04}
C.~P. Robert and G.~Casella, \emph{Monte {C}arlo {S}tatistical {M}ethods}, 2
  ed., Springer, 2004.

\bibitem[RCMA{\etalchar{+}}09]{ref:ramponinstable}
F.~Ramponi, D.~Chatterjee, A.~Milias-Argeitis, P.~Hokayem, and J.~Lygeros,
  \emph{Attaining mean square boundedness of a marginally stable noisy linear
  system with a bounded control input}, \url{http://arxiv.org/abs/0907.1436},
  2009.

\bibitem[RH05]{RichardsHow-05}
A.~Richards and J.~How, \emph{Robust model predictive control with imperfect
  information}, Proceedings of the American Control Conference, 2005,
  pp.~268--273.

\bibitem[SB09a]{SkafBoyd-Affine-09}
J.~Skaf and S.~Boyd, \emph{Design of affine controllers via convex
  optimization}, \url{http://www.stanford.edu/~boyd/papers/affine_contr.html},
  2009, To appear in IEEE Transactions on Automatic Control.

\bibitem[SB09b]{SkafBoyd-Qdesign-09}
\bysame, \emph{Nonlinear {Q}-design for convex stochastic control}, IEEE
  Transactions on Automatic Control \textbf{54} (2009), no.~10, 2426--2430.

\bibitem[SS85]{ref:schweitzerPolyApprox}
P.~J. Schweitzer and A.~Seidmann, \emph{Generalized polynomial approximations
  in {M}arkovian decision processes}, Journal of Mathematical Analysis and
  Applications \textbf{110} (1985), no.~2, 568--582.

\bibitem[SSW06]{ref:Stoorvogel}
A.~A. Stoorvogel, A.~Saberi, and S.~Weiland, \emph{On external semi-global
  stochastic stabilization of linear systems with input saturation},
  \url{http://homepage.mac.com/a.a.stoorvogel/subm03.pdf}, 2006.

\bibitem[ST03]{ref:tsiniasRobustStochISS}
J.~Spiliotis and J.~Tsinias, \emph{Notions of exponential robust stochastic
  stability, {ISS} and their {L}yapunov characterization}, International
  Journal of Robust and Nonlinear Control \textbf{13} (2003), no.~2, 173--187.

\bibitem[vHB03]{vanHessemFullSolution}
D.~H. van Hessem and O.~H. Bosgra, \emph{A full solution to the constrained
  stochastic closed-loop {MPC} problem via state and innovations feedback and
  its receding horizon implementation}, Proceedings of the 42nd IEEE Conference
  on Decision and Control, vol.~1, 2003, pp.~929--934.

\bibitem[vHB06]{vanHessem2006}
\bysame, \emph{Stochastic closed-loop model predictive control of continuous
  nonlinear chemical processes}, Journal of Process Control \textbf{16} (2006),
  no.~3, 225--241.

\bibitem[YB09]{ref:JunYanBitmead}
J.~Yan and R.~Bitmead, \emph{A constrained model-predictive approach to
  coordinated control}, To Appear in Automatica, 2009.

\bibitem[YSS97]{ref:YangSontagSussmann-97}
Y.~D. Yang, E.~D. Sontag, and H.~J. Sussmann, \emph{Global stabilization of
  linear discrete-time systems with bounded feedback}, Systems and Control
  Letters \textbf{30} (1997), no.~5, 273--281.

\end{thebibliography}

\newcommand{\etalchar}[1]{$^{#1}$}
\providecommand{\bysame}{\leavevmode\hbox to3em{\hrulefill}\thinspace}
\providecommand{\MR}{\relax\ifhmode\unskip\space\fi MR }
\providecommand{\MRhref}[2]{%
  \href{http://www.ams.org/mathscinet-getitem?mr=#1}{#2}
}
\providecommand{\href}[2]{#2}

\renewcommand{\theequation}{A.\arabic{equation}}
\setcounter{section}{1}  
\renewcommand{\thesection}{\Alph{section}}

\section*{Appendix}
	\subsection{Some identities}
	\label{s:facts}
	Recall the following standard special mathematical functions: the \emph{standard error function} $\erf(z) \Let \frac{2}{\sqrt\pi}\int_0^z \epower{-\frac{t^2}{2}}\mrm dt$ and the \emph{complementary error function}~\cite[p.~297]{ref:AbramowitzStegun} defined by $\erfc(z) \Let 1-\erf(z)$ for $z\in\R$, the \emph{incomplete Gamma function}~\cite[p.~260]{ref:AbramowitzStegun} defined by $\Gamma(a, z) \Let \int_z^\infty t^{a-1} \epower{-t} \mrm dt$ for $z, a > 0$, the \emph{confluent hypergeometric function}~\cite[p.~505]{ref:AbramowitzStegun} defined by $U(a, b, z) \Let \frac{1}{\Gamma(a)}\int_0^\infty \epower{-zt} t^{a-1} (1+t)^{b-a-1}\mrm dt$ for $a, b, z > 0$, and $\Gamma$ is the standard Gamma function. All of these are implemented as standard functions in Mathematica. The following facts can be found in \cite{ref:AbramowitzStegun} and are collected here for completeness.
	\begin{fact}
	\label{facts}
 		For $\sigma^2 > 0$ we have
		\begin{itemize}[leftmargin=*]
			\item $\displaystyle{\frac{1}{\sqrt{2\pi}\sigma}\int_z^\infty \epower{-\frac{t^2}{2\sigma^2}}\mrm dt = \frac{1}{2}\Bigl(1 + \erf\Bigl(\frac{z}{\sqrt 2\sigma}\Bigr)\Bigr)}$

			\item $\displaystyle{\frac{1}{\sqrt{2\pi}\sigma}\int_0^\infty \frac{t^2}{1+t^2} \epower{-\frac{t^2}{2\sigma^2}}\mrm dt = \frac{1}{2}\Bigl(\sqrt{2\pi}\sigma - \pi\epower{-\frac{1}{2\sigma^2}}\erfc\Bigl(\frac{1}{\sqrt 2\sigma}\Bigr)\Bigr)}$

			\item $\displaystyle{\frac{1}{\sqrt{2\pi}\sigma}\int_0^1 t^2 \epower{-\frac{t^2}{2\sigma^2}}\mrm dt = \sqrt{\frac{\pi}{2}}\sigma^3\erf\Bigl(\frac{1}{\sqrt 2\sigma}\Bigr) - \sigma^2\epower{-\frac{1}{2\sigma^2}}}$;

			\item $\displaystyle{\frac{1}{\sqrt{2\pi}\sigma}\int_1^\infty t \epower{-\frac{t^2}{2\sigma^2}}\mrm dt = \frac{\sigma}{\sqrt{2\pi}}\Gammaf(2\sigma^2, 1)}$

			\item $\displaystyle{\frac{1}{\sqrt{2\pi}\sigma}\int_0^\infty \frac{t^2}{\sqrt{1+t^2}} \epower{-\frac{t^2}{2\sigma^2}}\mrm dt = \frac{\sigma}{2\sqrt 2}U\Bigl(\frac{1}{2}, 0, \frac{1}{2\sigma^2}\Bigr)}$.
		\end{itemize}
	\end{fact}

\subsection{Proof of mean-square boundedness}
\label{s:app:proof}
%

\begin{proof}[Proof of Proposition \ref{prop:RHmain}]
	Fix $\xz\in\R^n$. For any $n\times n$ matrix $P = P\transp > 0$, using \eqref{eq:RHcl} and the fact that $\EE\left[\ee(w)\right]=0$, we see that for every $\ell = 1,\cdots, N_c$
	\begin{align*}
		\EE_{x_{kN_c}}\bigl[x_{kN_c+\ell}\transp P x_{kN_c + l}\bigr] & = x_{kN_c}\transp (A^\ell)\transp P A^\ell x_{kN_c}+2x_{kN_c}\transp (A^\ell)\transp PB_\ell \EE_{x_{kN_c}}\bigl[\pi^*_{kN_c:kN_c + \ell - 1}(x_{kN_c})\bigr]\\
		& \quad + \EE_{x_{kN_c}}\bigl[\norm{B_\ell\pi^*_{kN_c:kN_c + \ell - 1}(x_{kN_c}) + D_\ell\tilde w_{kN_c:kN_c + \ell - 1}}_P^2\bigr],
	\end{align*}
	where $\norm{\xi}_P \Let \sqrt{\xi\transp P\xi}$. Using the fact that $\norm{\pi^*_{kN_c:kN_c + \ell - 1}(x_{kN_c})}_\infty\le U_{\rm max}$ by construction, we obtain the following bound:
	\begin{align}
	\label{e:ellstep}
		\EE_{x_{kN_c}}\bigl[x_{kN_c+\ell}\transp P x_{kN_c+\ell}\bigr] \le  x_{kN_c}\transp (A^\ell)\transp P A^\ell x_{kN_c} +2 c_{1\ell} \norm{x_{kN_c}}_\infty + c_{2\ell},
	\end{align}
	where
	\begin{align*}
		c_{1\ell} & \Let m\norm{(A^\ell)\transp P B_\ell}_\infty U_{\rm max},\\
		c_{2\ell} & \Let m\norm{ B_\ell\transp P  B_\ell}_\infty U_{\max}^2 + \tr{D_\ell\transp P  D_\ell \Sigma_{ w}}\\
		& \quad +\max_{\tiny{\norm{\Upsilon(x_{kN_c})}_\infty\le U_{\max}/\phi_{\max}}}\big[\tr{ \Upsilon(x_{kN_c})\transp B_\ell\transp P B_\ell \Upsilon(x_{kN_c})\Lambda_1} +2\tr{\Upsilon(x_{kN_c})\transp B_\ell\transp P  D_\ell \Lambda_2}\big],\\
		& \text{and} \quad \Upsilon(x_{kN_c}) \Let \smat{\Theta^*_{1}(x_{kN_c})\\ \cdots\\\Theta^*_{N_c-1}(x_{kN_c})}.
	\end{align*}
	Since $A$ is a Schur stable matrix (and hence so is $A^\ell$) there exists \cite[Proposition 11.10.5]{ref:bernsteinMatrixMath} a matrix $P_\ell = P_\ell\transp > 0$ with real-valued entries that satisfies $(A^\ell)\transp P_\ell A^\ell - P_\ell = -\mathbf{I}_{n\times n}$; in particular, we have $x_{kN_c}\transp(A^\ell)\transp P_\ell A^\ell x_{kN_c} \le x_{kN_c}\transp P_\ell x_{kN_c} - x_{kN_c}\transp x_{kN_c}$. Therefore, with $P = P_\ell$ in \eqref{e:ellstep} we arrive at
	\begin{align}	
	\label{e:stab11}
		\EE_{x_{kN_c}}\bigl[x_{kN_c+\ell}\transp P_\ell x_{kN_c+\ell}\bigr] \le x_{kN_c}\transp P_\ell x_{kN_c} - \norm{x_{kN_c}}^2 + 2 c_{1\ell} \norm{x_{kN_c}}_\infty + c_{2\ell}.
	\end{align}
	For $\zeta_\ell\in\;]\max\{0,1-\lambda_{\max}(P_\ell)\}, 1[$ let $r_\ell \Let \frac{1}{\zeta_\ell}\bigl( c_{1\ell}+ \sqrt{c_{1\ell}^2 + c_{2\ell}\zeta_\ell}\bigr)$. Then elementary properties of the quadratic function $g(y) \Let -\zeta_\ell y^2 + 2c_{1\ell}y + c_{2\ell}$ show that
	\begin{align*}
		-\zeta_\ell \norm{x_{kN_c}}_\infty^2 + 2 c_{1\ell}\norm{x_{kN_c}}_\infty +  c_{2\ell} \le 0\quad\text{whenever }\norm{x_{kN_c}}_\infty >  r_\ell,
	\end{align*}
	In view of the above fact, simple manipulations in \eqref{e:stab11} lead to
	\begin{align*}
		\EE_{x_{kN_c}}\bigl[x_{kN_c+\ell}\transp P_\ell x_{kN_c + \ell}\bigr] \le x_{kN_c}\transp P_\ell x_{kN_c} - (1-\zeta_\ell)\norm{x_{kN_c}}^2 \quad \text{whenever } \norm{x_{kN_c}}_\infty > r_\ell,
	\end{align*}
	from which, letting $\rho_\ell \Let \Bigl(1- \frac{1-\zeta_\ell}{\lambdamax {P_\ell}}\Bigr)$, we arrive at
	\begin{align}
	\label{eqn:boundl}
		\EE_{x_{kN_c}}\bigl[x_{kN_c+\ell}\transp P_\ell x_{kN_c+\ell}\bigr] \le \rho_\ell x_{kN_c}\transp P_\ell x_{kN_c}\quad \text{whenever } \norm{x_{kN_c}}_\infty > r_\ell.
	\end{align}
	Let us define
	\begin{align*}
		\rho & \Let \max\limits_{\ell=1,\cdots,N_c}\rho_\ell, & r' &\Let\max\limits_{\ell=1,\cdots,N_c}r_\ell,\\
		\overline\lambda & \Let \max\limits_{\ell=1,\dots,N_c}\lambda_{\max}(P_\ell),  & \underline\lambda & \Let \min\limits_{\ell=1,\dots,N_c}\lambda_{\min}(P_\ell).
	\end{align*}
	Then we can obtain using~\eqref{eqn:boundl} the conservative bound for every $\ell = 1, \ldots, N_c$:
	\begin{equation*}
		\EE_{x_{kN_c}}\bigl[x_{kN_c+\ell}\transp P_{N_c} x_{kN_c+\ell}\bigr] \le \rho' x_{kN_c}\transp P_{N_c} x_{kN_c}\quad \text{whenever } \norm{x_{kN_c}}_\infty > r',
	\end{equation*}
	where $\rho'\Let\rho\frac{\overline\lambda\lambda_{\max}(P_{N_c})}{\underline\lambda \lambda_{\min}(P_{N_c})}$. It follows immediately that
	\begin{equation}
	\label{e:boundl2}
		\EE_{x_{kN_c}}\bigl[x_{kN_c+\ell}\transp P_{N_c} x_{kN_c+\ell}\bigr] \le \rho' x_{kN_c}\transp P_{N_c} x_{kN_c} + b'\indic{K'}(x_{kN_c}),
	\end{equation}
	where $K'\Let\bigl\{\xi\in\R^n\big|\norm{\xi}_\infty\le r'\bigr\}$.
	
	Let us define the function $V(\xi)\Let \xi\transp P_{N_c}\xi$, and fix $k\in\NN$ and $\ell =1,\dots,N_c$. Let $K_{N_c}\Let\bigl\{\xi\in\R^n\big|\norm{\xi}_\infty\le r_{N_c}\bigr\}$, $b\Let \sup\limits_{x\in K}\EE_x\bigl[V(x_{N_c})\bigr]$, and $b'\Let \max\limits_{\ell=1,\ldots, N_c}\sup\limits_{x\in K'}\EE_x\bigl[V(x_{\ell})\bigr]$. From \eqref{e:boundl2} we get
	\begin{align}
		\EE_{\xz}\bigl[V(x_{kN_c+\ell})\bigr] & = \EE_{\xz}\bigl[\EE\bigl[V(x_{kN_c+\ell})\,\big|\, x_{kN_c}\bigr]\bigr] \le \EE_{\xz}\bigl[\rho'V(x_{kN_c})+b'\indic{K'}(x_{kN_c})\bigr]\nn\\
		& \le \EE_{\xz}\bigl[\rho'\EE\bigl[V(x_{kN_c})\,\big|\, x_{(k-1)N_c}\bigr]+b'\indic{K'}(x_{kN_c})\bigr]\nn\\
		& \le \EE_{\xz}\bigl[\rho'\rho_{N_c} V(x_{(k-1)N_c})+b\indic{K_{N_c}}(x_{(k-1)N_c}) +b'\indic{K'}(x_{kN_c})\bigr]\nn\\
		& \cdots\nn\\
		& \le \rho'\rho_{N_c}^k V(x)+\sum_{i=0}^{k-1}b\rho_{N_c}^{k-1-i} \EE_{\xz}\bigl[\indic{K_{N_c}}(x_{iN_c})\bigr] +b'\EE_{\xz}\bigl[\indic{K'}(x_{kN_c})\bigr]\nn\\
		& \le \rho'\rho_{N_c}^k V(x)+\frac{b\bigl(1-\rho_{N_c}^k\bigr)}{1-\rho_{N_c}}+b'. \label{eqn:finalbound}
	\end{align}
	Note that the conditioning in the first few steps of~\eqref{eqn:finalbound} is well-defined because it is performed every $N_c$ steps starting from $0$, and the structure of our policy $\pi^*$ makes the process $(x_{tN_c})_{t\in\Nz}$ Markovian. Therefore, it follows from~\eqref{eqn:finalbound} that for all $t\Let kN_c+\ell$,
	\begin{align*}
		\sup\limits_{t\in\Nz}\EE_{\xz}\bigl[\norm{x_t}^2\bigr]& \le \frac{1}{\lambda_{\min}(P_{N_c})}\sup\limits_{t\in\Nz}\EE_{\xz}\bigl[V(x_{kN_c+\ell})\bigr]\\
		& \le \frac{1}{\lambda_{\min}(P_{N_c})}\left(\rho'\rho_{N_c}^k V(x)+\frac{b}{1-\rho_{N_c}}+b'\right)\\
		& <\infty,
	\end{align*}
	where the last step follows from the fact that $\rho_{N_c} < 1$. This completes the proof.
\end{proof}

\end{document}